\documentclass[11pt]{article}
\usepackage{graphicx}
\usepackage{amsmath,amsfonts,amssymb}
\usepackage{url}
\usepackage{amsmath}
\usepackage{graphicx}
\usepackage{epsfig}
\usepackage{epstopdf}
\usepackage{color}
\def\tto{\;{\lower 1pt \hbox{$\rightarrow$}}\kern -10pt
\hbox{\raise 2pt \hbox{$\rightarrow$}}\;}
\def\Hat{\widehat}

\def\Bar{\overline}
\def\ra{\rangle}
\def\la{\langle}
\def\ve{\varepsilon}
\def\epsilon{\varepsilon}
\def\B{\Bbb B}
\def\h{\hfill\Box}
\def\R{\Bbb R}

\def\N{I\!\!N}
\def\ox{\bar{x}}

\def\oy{\bar{y}}
\def\oz{\bar{z}}

\def\co{\mbox{\rm co}}

\def\graph{\mbox{\rm gph}\,}
\def\gph{\mbox{\rm gph}\,}

\def\epi{\mbox{\rm epi}\,}

\def\dom{\mbox{\rm dom}\,}

\def\i{\mbox{\rm int}\,}

\def\rge{\mbox{\rm rge}\,}
\def\i{\mbox{\rm int}\,}

\def\h{\hfill\square}

\def\O{\Omega}
\def\ph{\varphi}
\def\emp{\emptyset}

\def\oR{\Bar{\R}}
\def\lm{\lambda}
\def\gg{\gamma}
\def\dd{\delta}
\def\al{\alpha}
\def\bb{\beta}
\def\Th{\Theta}
\def\N{I\!\!N}

\newcounter{lk}

\begin{document}
\begin{center}
{\bf VARIATIONAL GEOMETRIC APPROACH TO GENERALIZED DIFFERENTIAL AND CONJUGATE CALCULI IN CONVEX ANALYSIS}\\[2ex]
B. S. MORDUKHOVICH\footnote{Department of Mathematics, Wayne State University, Detroit, MI 48202, USA (boris@math.wayne.edu) and RUDN University, Moscow 117198, Russia. Research of this author was partly supported by the National Science Foundation under grants DMS-1007132 and DMS-1512846 and by the Air Force Office of Scientific Research under grant \#15RT0462.}, N. M. NAM\footnote{Fariborz Maseeh Department of Mathematics and Statistics, Portland State University, Portland, OR 97207, USA (mau.nam.nguyen@pdx.edu). Research of this author was partly supported by the National Science Foundation under grant \#1411817.}, R. B. RECTOR\footnote{Fariborz Maseeh Department of Mathematics and Statistics, Portland State University, Portland, OR 97207, USA (r.b.rector@pdx.edu).} and T. TRAN\footnote{Fariborz Maseeh Department of Mathematics and Statistics, Portland State University, Portland, OR 97207, USA (tuyen2@pdx.edu).}.\\[2ex]
{\bf Dedicated to Michel Th\'era in honor of his 70th birthday}\vspace*{-0.1in}
\end{center}
\small{\bf Abstract.} This paper develops a geometric approach of variational analysis for the case of convex objects considered in locally convex topological spaces and also in Banach space settings. Besides deriving in this way new results of convex calculus, we present an overview of some known achievements with their unified and simplified proofs based on the developed geometric variational schemes.

\newtheorem{Theorem}{Theorem}[section]
\newtheorem{Proposition}[Theorem]{Proposition}
\newtheorem{Remark}[Theorem]{Remark}
\newtheorem{Lemma}[Theorem]{Lemma}
\newtheorem{Corollary}[Theorem]{Corollary}
\newtheorem{Definition}[Theorem]{Definition}
\newtheorem{Example}[Theorem]{Example}
\renewcommand{\theequation}{\thesection.\arabic{equation}}
\normalsize
\vspace*{-0.2in}

\section{Introduction}
\vspace*{-0.1in}

By now it has been widely accepted that convex analysis is one of the most important and useful areas of mathematical sciences with a great many applications
to optimization, economics, systems control, statistics, etc. The reader can find more information and extensive developments of various aspects of convex analysis and its applications in the books \cite{bc,bl,Ra,HU,kk,ph,r,r1,s,z} and the references therein. Michel Th\'era has made highly recognized contributions to the convex realm, particularly to the fields of convex duality, theory of monotone operators, and applications to optimization.\vspace*{-0.03in}

The major goal of this paper is to provide an adaptation and improvement of a certain {\em dual-space geometric} approach of {\em variational analysis} in the case of convex objects. Such an approach, based on the {\em extremal principle} for systems of closed (may not be convex) sets, has been systematically developed in \cite{m-book1} with broad applications to generalized differential calculus, optimization, equilibria, optimal control, economics, etc. However, several restrictions are imposed in the framework of \cite{m-book1} to proceed with the implementation of this approach. First, the local {\em closedness} of sets and {\em lower semicontinuity} (l.s.c.) of extended-real-valued functions, together with the {\em completeness} of the space in question,  are very essential in \cite{m-book1} as well as in other techniques of variational analysis; cf.\ \cite{BZ,RockWets-VA}. Furthermore, the main results of \cite{m-book1} hold in {\em Asplund} spaces, i.e., such Banach spaces where every separable subspace has a separable dual; see \cite{m-book1,ph} for more information on such spaces. \vspace*{-0.03in}

In this paper we develop a counterpart of the variational geometric approach to the study of convex sets, functions, and set-valued mappings in {\em locally convex topological vector} (LCTV) spaces without any completeness and closedness assumptions as well as in Banach spaces when the closedness is imposed. Elements of this approach to convex subdifferential calculus in finite-dimensional spaces have been presented in our previous publications \cite{bmn,bmn1} while some intersections rules in normed spaces have been very recently derived in \cite{mnam}. Here we address not only the {\em generalized differential calculus} but also calculus rules for the {\em Fenchel conjugate} of convex functions, which is a specific object of convex analysis.\vspace*{-0.03in}

The main variational geometric tool of our analysis in LCTV spaces is the {\em convex extremal principle} recently proposed in \cite{mnam} in the case normed spaces and extended in this paper to the general LCTV framework. Using it, we obtain refined calculus rules for supports and normals to set intersections and then apply this geometric machinery to developing other rules of generalized differential and conjugate calculi, while particularly concentrating on new issues to calculate subgradients and conjugates of the {\em marginal/optimal value} functions and (precise) {\em coderivative calculus} for convex-graph mappings. Besides this, in the case of Banach spaces under additional closedness assumptions, we exploit Attouch-Br\'ezis' technique \cite{AB} to deal with conjugate functions.\vspace*{-0.03in}

The rest of the paper is organized as follows. Section~2 recalls some definitions and preliminaries from convex and variational analysis. In Section~3 we study support functions to set intersections largely used in what follows. The key observation of this section is the extended version of the convex extremal principle applied to deriving a refined support and normal cone intersection rules in arbitrary LCTV spaces.\vspace*{-0.03in}

Section~4 applies the geometric approach together with the intersection results of Sections~3 to deriving major calculus rules for Fenchel conjugates of convex functions on LCTV and Banach spaces. In Section~5 we present new results on calculating in this way the Fenchel conjugate of the optimal value function in general convex settings.\vspace*{-0.03in}

Sections~6--8 are devoted to applying our variational geometric approach to generalized differential calculus for convex functions and set-valued mappings (multifunctions). We start in Section~6 with demonstrating the usage of this approach to give simple proofs of major sum and chain rules for subgradients and also some new developments for subdifferentiation of maximum functions. Section~8 presents new results on the exact subdifferential calculation for optimal value functions involving coderivatives of feasible sets. Finally, Section~9 develops new rules of coderivative calculus for convex-graph multifunctions.\vspace*{-0.03in}

The notation of this paper is standard in convex and variational analysis; cf.\ \cite{m-book1,r,RockWets-VA,z}. Throughout the paper we assume, unless otherwise stated, that all the spaces under consideration are separated LCTV. Given $\O\subset X$, we write $\R^+(\Omega):=\{tv\in X|\;t\in\R_+,\;v\in\Omega\}$, where $\R_+$ signifies the collection of positive numbers, and use the symbol $\overline\O$ for the topological closure of the set $\O$. The symbols $\B$ and $\B^*$ stand for the closed unit balls in the normed space $X$ in question and its topological dual $X^*$, respectively. In what follows we deal with {\em extended-real-valued} functions $f\colon X\to\oR:=(-\infty,\infty]$ and assume that they are {\em proper}, i.e., $\dom f:=\{x\in X|\;f(x)<\infty\}\ne\emp$.\vspace*{-0.2in}

\section{Basic Definitions and Preliminaries}
\setcounter{equation}{0}\vspace*{-0.1in}

In this section we briefly recall some basic notions for convex sets and functions studied in this paper, which can be found, e.g., in \cite{z} in the general LCTV space setting.\vspace*{-0.03in}

Given a convex set $\Omega\subset X$ with $\ox\in\Omega$, the {\em normal cone} to $\Omega$ at $\ox$ is
\begin{equation}\label{nor}
N(\ox;\Omega):=\big\{x^*\in X^*\big|\;\la x^*,x-\ox\ra\le 0\;\;\mbox{\rm for all }\;x\in\Omega\big\}.
\end{equation}
Given an extended-real-valued function $f\colon X\to\oR$ with its epigraph
$$
\epi f :=\big\{(x,\mu)\in X\times\R\big|\;\mu\ge f(x)\big\},
$$
we know that the standard notions of convexity and lower semicontinuity of $f$ can be equivalently described via the convexity and closedness of its epigraph, respectively. An element $x^*\in X^*$ is a {\em subgradient} of $f$ at $\ox\in\dom f$ if
\begin{equation}\label{sub}
\la x^*,x-\ox\ra\le f(x)-f(\ox)\;\;\mbox{\rm for all }\;x\in X,
\end{equation}
while the collection $\partial f(\ox)$ of all such $x^*$ is called the {\em subdifferential} of $f$ at $\ox$. Note that if $f\colon X^*\to\oR$, where $X^*$ is equipped with the weak$^*$ topology $\tau_{w^*}$, then $(X^*, \tau_{w^*})^*=X$. It follows that an element $x\in X$ is a subgradient of $f$ at $\ox^*\in X^*$ if and only if
\begin{equation*}
\la x^*-\ox^*,x\ra\leq f(x^*)-f(\ox^*)\; \; \mbox{\rm for all }x^*\in X^*.
\end{equation*}
It is easy to check that the subdifferential $\partial f(\ox)$ can be represented geometrically as
\begin{equation*}
\partial f(\ox)=\big\{x^*\in X^*\big|\;(x^*,-1)\in N\big((\ox,f(\ox));\epi f\big)\big\}.
\end{equation*}
On the other hand, we have $N(\ox;\O)=\partial\dd(\ox;\O)$, where $\dd(x;\O)=\dd_\O(x)$ is the indicator function of the set $\O$ equal to $0$ if $x\in\O$ and to $\infty$ if $x\notin\O$.\vspace*{-0.03in}

The {\em Fenchel conjugate} $f^*\colon X^*\to[-\infty,\infty]$ of $f\colon X\to\oR$ is defined by
\begin{equation}\label{conj}
f^*(x^*):=\sup\big\{\la x^*,x\ra-f(x)\big|\;x\in X\big\}.
\end{equation}
Note that if $f$ is proper, then $f^*\colon X^*\to \oR$ is convex and l.s.c.\ on $X^*$ regardless of the convexity and lower semicontinuity of $f$. Recall also the {\em second Fenchel conjugate} (or {\em biconjugate}) of $f\colon X\to\oR$ given by
\begin{equation*}
f^{**}(x):=(f^*)^*(x)=\sup\big\{\la x^*,x\ra-f^*(x^*)\big|\;x^*\in X^*\big\},\quad x\in X.
\end{equation*}
As well known, $f^{**}=f$ on $X$ provided that $f$ is proper, l.s.c., and convex.\vspace*{-0.03in}

Another important operation on functions is the infimal convolution. Given two functions $f,g: X\to\oR$, their {\em infimal convolution} $f\oplus g$ is defined by
\begin{equation}\label{inf-con}
(f\oplus g)(x):=\inf\big\{f(x_1)+g(x_2)\big|\; x_1+x_2=x\big\}=\inf\big\{f(u)+g(x-u)\big|\;u\in X\big\},\quad x\in X.
\end{equation}
Assuming that $(f\oplus g)(x)>-\infty$ for all $x\in X$, we clearly have the convexity of the infimal convolution $f\oplus g$ provided that both functions $f,g$ are convex.\vspace*{-0.03in}

Finally in this section, recall a generalized differentiation notion for set-valued mappings, which has not been investigated in the standard framework of convex analysis while has been well recognized in the general setting of variational analysis \cite{m-book1,RockWets-VA}. Given a set-valued mapping $F\colon X\tto Y$  between LCTV spaces with its domain and graph
$$
\dom F:=\big\{x\in X\big|\;F(x)\ne\emp\big\}\;\;\mbox{\rm and }\;\gph F:=\big\{(x,y)\in X\times Y\big|\;y\in F(x)\big\},
$$
we define the {\em coderivative} $D^*F(\ox,\oy)\colon Y^*\tto X^*$ of $F$ at $(\ox,\oy)\in\gph F$ as a positively homogeneous multifunction on $Y^*$ taken the values
\begin{equation}\label{cod}
D^*F(\ox,\oy)(y^*):=\big\{x^*\in X^*\big|\;(x^*,-y^*)\in N\big((\ox,\oy);\graph F\big)\big\},\quad y^*\in Y^*.
\end{equation}
If $F(x):=Ax$ is a continuous linear operator $A\colon X\to Y$ with $\oy:=A\ox$, then
\begin{equation*}
D^*A(\ox,\oy)(y^*)=\big\{A^*y^*\big\}\;\;\mbox{\rm for all }\;\ox\in X,\;\;y^*\in Y^*,
\end{equation*}
via the adjoint operator $A^*\colon Y^*\to X^*$ given by $A^*y^*=y^*\circ A$ for $y^*\in Y^*$. This means that the coderivative is a set-valued extension of the classical notion of the adjoint linear operator in functional analysis.\vspace*{-0.2in}

\section{Extremal Principle and Supports to Set Intersections}
\setcounter{equation}{0}\vspace*{-0.1in}

In this section we formulate the enhanced notion of set extremality, which is the most appropriate for convex sets in LCTV spaces, and then establish the {\em convex extremal principle} in this setting. The main applications of the convex extremal principle are given here to deriving refined support and normal cone intersection rules.

We start with some auxiliary properties of support functions needed in what follows. Recall that the {\em support function} $\sigma_\O\colon X^*\to\oR$ of a nonempty subset $\O\subset X$  is defined by
\begin{equation}\label{supp}
\sigma_\O(x^*):=\sup\big\{\la x^*,x\ra\big|\;x\in\O\big\},\quad x^*\in X^*.
\end{equation}
Note that $\sigma_\O$ is always convex on $X^*$ regardless of the convexity of $\O$ and that we have the conjugacy relationship $\sigma_\O(x^*)=\dd^*_\O(x^*)$ as $x^*\in X^*$ with the indicator function of $\O$.\vspace*{-0.03in}

Let us first present the following known result widely used below. We provide a direct proof of it for the reader's convenience.

\begin{Proposition}{\bf(conjugates to support function convolutions via set intersections).}\label{p1} Let the sets $\O_1$ and $\O_2$ be nonempty, closed, and convex in $X$ with $\O_1\cap\O_2\ne\emp$. Then we have the representation
\begin{equation}\label{i1}
\big(\sigma_{\O_1}\oplus\sigma_{\O_2}\big)^*(x)=\delta_{\O_1\cap\O_2}(x)\;\;\mbox{\rm for all }\;x\in X.
\end{equation}
\end{Proposition}\vspace*{-0.12in}
{\bf Proof.} Consider the infimal convolution \eqref{inf-con} of the support functions of $\O_i$ given by
\begin{equation}\label{inf-supp}
f(x^*):=\big(\sigma_{\O_1}\oplus\sigma_{\O_2}\big)(x^*),\quad x^*\in X^*.
\end{equation}
It follows from \cite[Corollay~2.4.7]{r} and a well-known subdifferential result for support functions  on LCTV spaces that
\begin{equation*}
\partial f(0)=\partial\sigma_{\O_1}(0)\cap\partial\sigma_{\O_2}(0)=\O_1\cap\O_2.
\end{equation*}
It is easy to see that $f(x^*)$ is convex and positively homogeneous with $f(0)=0$ under the assumptions made, and so by \cite[Theorem~2.4.11(i)]{z} we get
\begin{equation*}
f^*(x)=\dd_{\partial f(0)}(x)\;\;\mbox{\rm whenever }\;x\in X,
\end{equation*}
which readily verifies the infimal convolution representation in \eqref{i1}.$\h$

The proofs of the next two lemmas follow the arguments from \cite{AB} with some simplifications in the case of the support functions under consideration. They are essentially based on the imposed Banach space and closedness assumptions.\vspace*{-0.12in}

\begin{Lemma}{\bf(weak$^*$ compactness of support level sets).}\label{lm1} Let $\O_1$ and $\O_2$ be nonempty closed subsets of a Banach space $X$. Suppose that $\R^+(\O_1-\O_2)=X$. Then for any numbers $\al,\bb\ge 0$ the restricted level set
\begin{equation*}
K=K_{\alpha,\beta}:=\big\{(x^*_1,x^*_2)\in X^*\times X^*\big|\;\sigma_{\O_1}(x^*_1)+\sigma_{\O_2}(x^*_2)\le\al,\;\|x^*_1+x^*_2\|\le\bb\big\}
\end{equation*}
is compact in the weak$^*$ topology of the dual product space $X^*\times X^*$.
\end{Lemma}\vspace*{-0.12in}
{\bf Proof.} The closedness of the set $K$ in the weak$^*$ topology of $X^*\times X^*$ is obvious. By the classical Alaoglu-Bourbaki theorem it remains to show that this set is norm-bounded in $X^*\times X^*$. Having in mind the uniform boundedness principle, we need to verify that the set of linear continuous functionals from $K$ is bounded pointwise. To proceed, take any $(x_1,x_2)\in X\times X$ and find by the imposed assumption $\R^+(\O_1-\O_2)=X$ elements $\lambda>0$, $w_1\in\O_1$, and $w_2\in\O_2$ for which $x_1-x_2=\lambda(w_1-w_2)$. Then we have
\begin{align*}
\la x^*_1,x_1\ra+\la x^*_2,x_2\ra&=\lambda\la x^*_1,w_1\ra+\lambda\la x^*_2,w_2\ra+\la x^*_1+x^*_2,x_2-\lambda w_2\ra\\
&\le\lambda\big(\sigma_{\O_1}(x^*_1)+\sigma_{\O_2}(x^*_2)\big)+\|x^*_1+x^*_2\|\cdot\|x_2-\lambda w_2\|\le\lambda\al+\bb\|x_2-\lambda w_2\|.
\end{align*}
Since the latter also holds for $(-x_1,-x_2)$, we arrive at the conclusion of the lemma. $\h$\vspace*{-0.12in}

\begin{Lemma}{\bf(lower semicontinuity of support infimal convolutions).}\label{lm2} In the setting of Lemma~{\rm\ref{lm1}} we have that the infimal convolution
$(\sigma_{\O_1}\oplus\sigma_{\O_2})\colon X^*\to\oR$ is lower semicontinuous with respect to the weak$^*$ topology of $X^*$.
\end{Lemma}\vspace*{-0.12in}
{\bf Proof.} It suffices to prove that for any $\gamma\in\R$ the set
\begin{equation*}
C:=\big\{x^*\in X^*\big|\;\big(\sigma_{\O_1}\oplus\sigma_{\O_2}\big)(x^*)\le\gamma\big\}
\end{equation*}
is weak$^*$ closed in $X^*$. Consider the parametric set family
\begin{equation*}
C_\epsilon:=\big\{x^*\in X^*\big|\;x^*=x^*_1+x^*_2,\;\sigma_{\O_1}(x^*_1)+\sigma_{\O_2}(x^*_2)\le\gamma+\epsilon\big\},\quad\ve>0,
\end{equation*}
with $C=\bigcap_{\epsilon>0}C_\epsilon$ and show that each set $C_\epsilon$ as $\epsilon>0$ is weak$^*$ closed in $X^*$. Using the seminal Banach-Dieudonn\'e-Krein-\u Smulian theorem (see, e.g., \cite[Theorem~V.5.7]{ds}), we only need to check that the intersection $C_{\epsilon}\cap r\B^*$ is weak$^*$ closed in $X^*$ for all $r>0$. To this end, define the mapping $T\colon X^*\times X^*\to X^*$ by
\begin{equation*}
T(x^*_1,x^*_2)=x^*_1+x^*_2,
\end{equation*}
which is continuous in the weak$^*$ topology of $X^*\times X^*$. It is easy to observe that
\begin{equation*}
C_{\epsilon}\cap r\B^*=T(K_{\gamma+\epsilon,r}),
\end{equation*}
where $K_{\gamma+\epsilon,r}$ is defined in Lemma~\ref{lm1}. It follows from Lemma~\ref{lm1} that $C_{\epsilon}\cap r\B^*$ is weak$^*$ compact and hence weak$^*$ closed, which verifies that $\sigma_{\O_1}\oplus\sigma_{\O_2}$ is weak$^*$ l.s.c.\ on $X^*$. $\h$

To proceed further, let us first formulate an appropriate LCTV version of set extremality defined recently in \cite{mnam} in normed spaces. It seems to be more efficient to deal with convex sets in LCTV spaces in comparison with its local extremality counterpart developed in \cite{m-book1}.\vspace*{-0.12in}

\begin{Definition}{\bf(set extremality in LCTV spaces).}\label{extrem} We say that two nonempty sets $\O_1,\O_2\subset X$ form an {\sc extremal system} in the space $X$ if for any neighborhood $V$ of the origin there exists a vector $a\in X$ such that
\begin{equation}\label{ext}
a\in V\;\;\mbox{\rm and }\;(\O_1+a)\cap\O_2=\emp.
\end{equation}
\end{Definition}

Note that \eqref{ext} does note require that the sets $\O_1,\O_2$ have the nonempty intersection in contrast to \cite{m-book1}. As discussed in \cite{mnam} similarly to \cite[Section~2.1]{m-book1} for the local counterpart, the extremality notion \eqref{ext} covers various concepts of optimal solutions to problems of scalar and multiobjective optimization, equilibrium systems, etc. with numerous applications.\vspace*{-0.02in}

The following result is an extension of \cite[Theorem~2.2]{mnam} for the case of normed spaces.\vspace*{-0.12in}

\begin{Theorem}{\bf(convex extremal principle).}\label{cep} Let the sets $\O_1,\O_2\subset X$ be nonempty and convex. Then we have the assertions: \vspace*{-0.03in}

{\bf(i)} The sets $\O_1$ and $\O_2$ form an extremal system in $X$ if and only if $0\notin{\rm int}(\Omega_1-\Omega_2)$. Furthermore, the extremality of $\O_1,\O_2$ implies that $({\rm int}\,\O_1)\cap\O_2=\emp$ and $({\rm int}\,\O_2)\cap\O_1=\emp$.\vspace*{-0.03in}

{\bf(ii)} If $\O_1$ and $\O_2$ form an extremal system in $X$ and the set difference $\O_1-\O_2$ is solid, i.e.,
\begin{equation}\label{solid}
{\rm int}(\O_1-\O_2)\ne\emp,
\end{equation}
then the set $\O_1$ and $\O_2$ are separated, i.e.,
\begin{equation}\label{sep}
\sup_{x\in\O_1}\la x^*,x\ra\le\inf_{x\in\O_2}\la x^*,x\ra\;\;\mbox{\rm for some }\;x^*\in X^*,\;x^*\ne 0.
\end{equation}\vspace*{-0.2in}

{\bf(iii)} The separation property \eqref{sep} always yields the set extremality \eqref{ext} even if \eqref{solid} fails.
\end{Theorem}\vspace*{-0.12in}
{\bf Proof.} {\rm\bf (i)} Let the sets $\O_1$ and $\O_2$ form an extremal system in $X$. Arguing by contradiction, suppose that $0\in{\rm int}(\Omega_1-\Omega_2)$. Then there is a balanced neighborhood $V$ of $0\in X$ with
\begin{equation*}
V\subset\O_1-\O_2.
\end{equation*}
For any $a\in V$ we have $-a\in V\subset\O_1-\O_2$ and arrive at the contradiction $(\O_1+a)\cap\O_2\ne\emp$.

Now suppose that $0\notin{\rm int}(\Omega_1-\Omega_2)$. Then for any neighborhood $V$ of the origin we get
$$
V\cap\big[X\setminus(\O_1-\O_2)\big]\ne\emp.
$$
Assume without loss of generality that $V$ is balanced, which yields $(-V)\cap[X\setminus(\O_1-\O_2)]\ne\emp$. Choosing $a\in V$ gives us the inclusion
$$
-a\in\big[X\setminus(\O_1-\O_2)\big],
$$
and hence $(\O_1+a)\cap\O_2=\emp$. This verifies that the sets $\O_1$ and $\O_2$ form an extremal system.

For the second part of (i), suppose by contradiction that $({\rm int}\,\O_1)\cap\O_2\ne\emp$. Then there exists a vector $x\in{\rm int}\,\O_1$ with $x\in\O_2$. We can always choose a balanced neighborhood $V$ of the origin such that $x+V\subset\O_1$. For any $a\in V$ it shows that $-a\in V$ and $x-a\in\O_1$. Hence $(a+\O_1)\cap\O_2\ne\emp$, which is a contradiction.

{\rm\bf (ii)} If $\O_1$ and $\O_2$ form an extremal system, then it follows from (i) that $0\notin{\rm int}(\Omega_1-\Omega_2)$. In addition to this, the assumption of (ii) on the solidness of the set difference $\O_1-\O_2$ allows us to use the convex separation theorem, which yields \eqref{sep}.

{\rm\bf (iii)} Suppose that \eqref{sep} holds, which gives us a vector $c\in X$ with $\la x^*,c\ra>0$. Fix any neighborhood $V$ of the origin. Since $V$ is always absorbing, we can select a natural number $k$ sufficiently large such that $a:=-c/k\in V$. Let us show that \eqref{ext} is satisfied with this vector $a$. Indeed, the negation of this means that then there exists $\Hat x\in\O_2$ such that $\Hat x-a\in\O_1$. By the separation property from \eqref{sep} we have
$$
\la x^*,\Hat x-a\ra\le\sup_{x\in\O_1}\la x^*,x\ra\le\inf_{x\in\O_2}\la x^*,x\ra\le\la x^*,\Hat x\ra,
$$
On the other hand, the construction of $a$ tells us that
$$
\la x^*,\Hat x\ra-\la x^*,a\ra=\la x^*,\Hat x\ra+k\la x^*,c\ra\le\la x^*,\Hat x\ra\;\mbox{ for all }\;k\in\N.
$$
This shows that $\la x^*,c\ra\le 0$, which contradicts the above choice of the vector $c$. $\h$

It is worth mentioning that \cite[Theorem~2.2]{mnam} offers a version of Theorem~\ref{cep}(ii) in Banach spaces with the replacement of \eqref{solid} by the {\em sequential normal compactness} (SNC) property imposed on one of the sets $\O_1,\O_2$. The SNC property is effectively characterized in \cite[Theorem~1.21]{m-book1} for convex sets. Furthermore, in \cite[Theorem~2.5]{mnam} the reader can find an {\em approximate} version of the convex extremal principle without imposing either \eqref{cep} or the SNC property. The latter approximate version yields, in particular, the celebrated Bishop-Phelps theorem establishing the density of support points on the boundary of any closed convex subset of a Banach space; see \cite[Theorem~3.18]{ph}.\vspace*{-0.03in}

The next theorem presents major geometric results on representing the support function for convex set intersections via the infimal convolution of the support functions for each component. We derive these representations under three generally independent groups of assumptions, which constitute the corresponding assertions of the theorem. The first assertion is established in the general LCTV setting without any closedness requirements on the sets in question under the {\em difference interiority condition}
\begin{equation}\label{dqc}
0\in{\rm int}(\O_1-\O_2)
\end{equation}
and the boundedness of one of the sets by using the convex extremal principle from Theorem~\ref{cep}. The second assertion, which is proved similarly in any LCTV space, does not impose any set boundedness while requires a more restrictive qualification condition; see more discussions after the proof of the theorem. The third assertion holds for closed sets in Banach spaces with replacing \eqref{dqc} by the weaker {\em Attouch-Br\'ezis qualification condition}
\begin{equation}\label{ab}
\R^+(\O_1-\O_2)\;\mbox{ is a closed subspace of }\;X
\end{equation}
by using the above lemmas and Proposition~\ref{p1}.\vspace*{-0.03in}

Recall that a subset $\O$ of an LCTV space $X$ is {\em bounded} if for any neighborhood $V$ of the origin there exists $\gamma>0$ such that $\O\subset \alpha V$ whenever $|\alpha|>\gamma$.\vspace*{-0.12in}

\begin{Theorem}{\bf(support functions to set intersections).}\label{sigma intersection rule} Let the sets $\Omega_1,\Omega_2\subset X$ be nonempty and convex.
Suppose that   one of the following conditions is satisfied:\\[1ex]
{\bf(i)} The difference interiority condition \eqref{dqc} is satisfied and the set $\O_2$ is bounded.\\
{\bf (ii)} Either $({\rm int}\,\O_2)\cap\O_1\neq\emp$, or $({\rm int}\,\O_1)\cap\O_2\ne\emp$.\\
{\bf(iii)} The space $X$ is Banach, both sets $\Omega_1$ and $\Omega_2$ are closed, and the Attouch-Br\'ezis qualification condition \eqref{ab} holds.\\[1ex]
Then the support function \eqref{supp} to the intersection $\O_1\cap\O_2$ is represented as
\begin{equation}\label{supp1}
\big(\sigma_{\Omega_1\cap\Omega_2}\big)(x^*)=\big(\sigma_{\Omega_1}\oplus\sigma_{\Omega_2}\big)(x^*)\;\;\mbox{\rm for all }\;x^*\in X^*.
\end{equation}
Furthermore, for any $x^*\in\dom(\sigma_{\Omega_1\cap\Omega_2})$ there are $x^*_1,x^*_2\in X^*$ such that $x^*=x^*_1+x^*_2$ and
\begin{equation}\label{supp2}
(\sigma_{\Omega_1\cap\Omega_2})(x^*)=\sigma_{\Omega_1}(x^*_1)+\sigma_{\Omega_2}(x^*_2).
\end{equation}
\end{Theorem}\vspace*{-0.12in}
{\bf Proof.} First we justify both formula \eqref{supp1} and representation \eqref{supp2} in case (i). To verify the inequality ``$\le$" in \eqref{supp1},
fix $x^*\in X^*$ and pick $x^*_1,x^*_2$ with $x^*=x^*_1+x^*_2$. Then
\begin{align*}
\la x^*,x\ra=\la x^*_1,x\ra+\la x^*_2,x\ra\le\sigma_{\O_1}(x^*_1)+\sigma_{\O_2}(x^*_2)\;\;\mbox{\rm for any }\;x\in\O_1\cap\O_2.
\end{align*}
Taking the infimum on the right-hand side above with respect to all such $x^*_1,x^*_2$ yields
\begin{equation*}
\la x^*,x\ra\le\big(\sigma_{\O_1}\oplus\sigma_{\O_2}\big)(x^*)\;\;\mbox{\rm whenever }\;x\in\O_1\cap\O_2,
\end{equation*}
which clearly implies the inequality ``$\le$" in \eqref{supp1} without imposing the assumptions in (i).\vspace*{-0.03in}

Next we prove the inequality ``$\ge$" in \eqref{supp1} under the assumptions in (i). To proceed, fix any $x^*\in\dom(\sigma_{\Omega_1\cap\Omega_2})$, denote $\alpha:=\sigma_{\O_1\cap\O_2}(x^*)\in\R$ for which we have
\begin{equation*}
\la x^*,x\ra\le\alpha\;\;\mbox{\rm whenever }\;x\in\O_1\cap\O_2,
\end{equation*}
and define the nonempty convex subsets of $X\times\R$ by
\begin{equation}\label{theta}
\Theta_1:=\O_1\times[0,\infty)\;\;\mbox{\rm and }\;\Theta_2:=\big\{(x,\lambda)\in X\times\R\big|\;x\in\O_2,\;\lambda\le\la x^*,x\ra-\alpha\big\}.
\end{equation}
It is easy to see from the constructions of $\Th_1$ and $\Th_2$ that
\begin{equation*}
\big(\Theta_1+(0,\gamma)\big)\cap\Theta_2=\emp\;\;\mbox{\rm for any }\;\gg>0,
\end{equation*}
and thus these sets form and {\em extremal system} in $X\times\R$. We deduce from Theorem~\ref{cep}(i) that $0\notin\mbox{\rm int}(\Theta_1-\Theta_2)$. To apply Theorem~\ref{cep}(ii) to these sets, we need to verify that the set difference $\Theta_1-\Theta_2$ is solid. The property $0\in\mbox{\rm int}(\O_1-\O_2)$ allows us to a neighborhood $U$ of the origin such that $U\subset \O_1-\O_2$. By the continuity of the function $\ph(x):=\la -x^*, x\ra +\alpha$ for $x\in X$, there exists a neighborhood $W$ of the origin such that $\ph$ is bounded above on $W$. Since $\O_2$ is bounded, we can find $t>0$ such that $\O_2\subset tW$. Note that $\ph$ is also bounded above on $tW$, so we can find $\bar{\lambda}\in \R$ such that
\begin{equation*}\label{lambda}
\bar{\lambda}\ge\sup_{x\in tW}\la-x^*,x\ra+\alpha\ge\sup_{x\in \O_2}\la-x^*,x\ra+\alpha.
\end{equation*}
Let us check that $U\times(\bar{\lambda},\infty)\subset\Theta_1-\Theta_2$, and so \eqref{solid} holds. Indeed, for $(x,\lambda)\in U\times(\bar{\lambda},\infty)$ we have $x\in U\subset \O_1-\O_2$ and $\lambda>\bar{\lambda}$. Hence $x=w_1-w_2$ with some $w_1\in\O_1$, $w_2\in\O_2$ and
\begin{equation*}
(x,\lambda)=(w_1,\lambda-\bar\lambda)-(w_2,-\bar\lambda).
\end{equation*}
It follows from $\lambda-\bar\lambda>0$ that $(w_1,\lambda-\bar\lambda)\in\Theta_1$. We deduce from \eqref{theta} and the choice of $\bar\lm$ that $(w_2,-\bar\lambda)\in\Theta_2$, and thus $\mbox{\rm int}(\Theta_1-\Theta_2)\neq\emp$. Theorem~\ref{cep}(ii) tells us now that there exists a pair $(0,0)\ne(y^*,\beta)\in X^*\times\mathbb{R}$ for which we get
\begin{equation}\label{sep f}
\la y^*,x\ra+\beta\lambda_1\le\la y^*,y\ra+\beta\lambda_2\;\;\mbox{\rm whenever }\;(x,\lambda_1)\in\Theta_1,\;(y,\lambda_2)\in\Theta_2.
\end{equation}
It follows from the structure of $\Th_1$ that $\bb\le 0$. Assuming that $\bb=0$ tells us by \eqref{sep f} that
\begin{equation*}
\la y^*,x\ra\le\la y^*,y\ra\;\;\mbox{\rm for all }\;x\in\Omega_1,\;y\in\Omega_2.
\end{equation*}
This yields $y^*=0$ due to $0\in {\rm int}(\Omega_1-\Omega_2)$, a contradiction. Take now the pairs $(x,0)\in\Theta_1$ and $(y,\langle x^*,y\rangle-\alpha)\in\Theta_2$ in \eqref{sep f} and get
\begin{equation*}
\la y^*,x\ra\le\la y^*,y\ra+\beta(\la x^*,y\ra-\alpha)\;\;\mbox{\rm with }\;\bb<0,
\end{equation*}
which brings us to the estimate
\begin{equation*}
\alpha\ge\big\la y^*/\beta+x^*,y\big\ra+\big\la-y^*/\bb,x\big\ra\;\;\mbox{\rm for all }\;x\in\Omega_1,\;y\in\Omega_2.
\end{equation*}
By putting $x^*_1:=y^*/\beta+x^*$ and $x^*_2:=-y^*/\beta$ we arrive at the inequality ``$\ge$" in \eqref{supp1} and representation \eqref{supp2}. This justifies the conclusions of the theorem in case (i).\vspace*{-0.03in}

To verify the results under (ii), it clearly suffices to examine only the first case therein. Considering the sets $\Th_1,\Th_2$ from \eqref{theta}, we see that
\begin{equation*}
{\rm int}\,\Th_2=\big\{(x,\lm)\in X\times\R\big|\;x\in\i\O_2,\;\lm<\la x^*,x\ra-\al\big\}\ne\emp,
\end{equation*}
and so ${\rm int}(\Th_1-\Th_2)\ne\emp$. Furthermore, it follows from the assumption $(\i\O_2)\cap\O_1\ne\emp$ that $0\in{\rm int}(\O_1-\O_2)$, which allows us to proceed as in the above proof in case (i).\vspace*{-0.03in}

Consider finally the Banach space case (iii) of the theorem and first verify its conclusions when $\R^+(\O_1-\O_2)=X$. Taking the Fenchel conjugate in both sides
of formula \eqref{i1} from Proposition~\ref{p1}, which holds even without the assumptions of (ii),
and then using Lemma~\ref{lm2} gives us the equalities
\begin{equation*}
\dd^*_{\O_1\cap\O_2}(x^*)=\sigma_{\O_1\cap\O_2}(x^*)=\big(\sigma_{\O_1}\oplus\sigma_{\O_2}\big)^{**}(x^*)=\big(\sigma_{\O_1}\oplus\sigma_{\O_2}\big)(x^*)
\end{equation*}
for all $x^*\in X^*$, which justify the representation in \eqref{supp1} when the assumption $\R^+(\O_1-\O_2)=X$ is satisfied. Treating further the general case of \eqref{ab} in (iii), denote $L:=\R^+(\O_1-\O_2)$ the closed subspace of $X$ in question. Since $\O_1\cap\O_2\ne\emp$ by \eqref{ab}, we can always translate the situation to $0\in\O_1\cap\O_2$, and hence suppose that $\O_1,\O_2\subset L$. This reduces the general case under \eqref{ab} to the one $\R^+(\O_1-\O_2)$ treated above; so \eqref{supp1} is justified.\vspace*{-0.03in}

Representation \eqref{supp2} under (iii) for $x^*\in\dom(\sigma_{\O_1}\cap\sigma_{\O_2})$ follows from the weak$^*$ compactness of the set $K_{\al,\bb}$ in Lemma~\ref{lm1} with $\al:=(\sigma_{\O_1}\oplus\sigma_{\O_2})(x^*)+\epsilon$, where $\epsilon>0$ is arbitrary and where $\bb:=\|x^*\|$. This completes the proof of the theorem. $\h$

The results of type \eqref{supp2} in Theorem~\ref{sigma intersection rule} (i.e., ensuring that the infimal convolution is {\em exact}) go back to Moreau \cite{more} in Hilbert spaces. More recently \cite{Ernst-Thera07} Ernst and Th\'era established  necessary and sufficient conditions for this property in finite-dimensional spaces.\vspace*{-0.03in}

If $X$ is Banach and both sets $\O_i$ are closed with $\mbox{\rm int}(\O_1-\O_2)\ne\emp$, then \eqref{dqc} reduces to the {\em core qualification condition} $0\in{\rm core}(\O_1-\O_2)$ developed by Rockafellar \cite{r1}, where
\begin{equation*}
{\rm core}\,\Omega:=\big\{x\in\Omega\big|\;\forall\,v\in X\;\exists\,\gg>0\;\mbox{\rm such that }\;x+tv\in\Omega\;\mbox{\rm whenever }\;|t|<\gg\big\}.
\end{equation*}
This follows from the well-known facts that $\i\O$=${\rm core}\,\O$ for closed convex subsets of Banach spaces (see, e.g., \cite[Theorem~4.1.8]{BZ}) and that $\i\,\Bar\O$=$\i\,\O$ provided that $\i\,\O\ne\emp$. Note that the Attouch-Br\'ezis condition \eqref{ab} from \cite{AB} essentially supersedes the equivalent core and difference interior qualification conditions for closed sets in Banach spaces, while the latter works in general LCTV spaces. As shown in \cite{mnam}, the conventional interiority conditions in (ii) of Theorem~\ref{sigma intersection rule} {\em strictly} imply the conditions in (i), including the set boundedness, provided that the space $X$ is {\em normed}. We cannot conclude this in general, unless the LCTV space $X$ in question has a basis of {\em bounded} neighborhoods of the origin.\vspace*{-0.03in}

The next result provides various qualification conditions in LCTV and Banach space settings under which the important normal cone intersection formula holds. This result under the qualification condition \eqref{qc} has been derived in \cite{mnam} directly from the convex extremal principle in normed spaces. Below we present a unified derivation of the normal cone intersection rule from the corresponding conditions of Theorem~\ref{sigma intersection rule}.\vspace*{-0.12in}

\begin{Theorem}{\bf (normal cone intersection rule).}\label{nir} Let the sets $\Omega_1,\Omega_2\subset X$ be convex with $\ox\in\Omega_1\cap\Omega_2$. Suppose that one of the following conditions {\rm(i)}--{\rm(iii)} is satisfied:\\[1ex]
{\bf(i)} There exists a bounded convex neighborhood $V$ of $\ox$ such that
\begin{equation}\label{qc}
0\in\mbox{\rm int}\big(\Omega_1-(\Omega_2\cap V)\big).
\end{equation}
{\bf (ii)} Either $({\rm int}\,\O_2)\cap\O_1\ne\emp$, or $({\rm int}\,\O_1)\cap\O_2\ne\emp$.\\
{\bf(iii)} The space $X$ is Banach, both sets $\Omega_1$ and $\Omega_2$ are closed, and the Attouch-Br\'ezis qualification condition \eqref{ab} holds.\\[1ex]
Then we have the normal cone intersection rule
\begin{equation}\label{ni}
N(\ox;\Omega_1\cap\Omega_2)=N(\ox;\Omega_1)+N(\ox;\Omega_2).
\end{equation}
\end{Theorem}\vspace{-0.1in}
{\bf Proof.} First we verify \eqref{ni} under (i). Denote $A:=\Omega_1$, $B:=\Omega_2\cap V$ and observe by \eqref{qc} that $0\in\mbox{\rm int}(A-B)$ and $B$ is bounded. It follows from the normal cone definition \eqref{nor} that $x^*\in N(\ox;\O)$ for $\ox\in\O$ if and only if $\sigma_{\O}(x^*)=\la x^*,\ox\ra$. Then pick $x^*\in N(\ox;A\cap B)$ and get $\la x^*,\ox\ra=\sigma_{A\cap B}(x^*)$. By Theorem~\ref{sigma intersection rule}(i) there are $x^*_1,x^*_2\in X^*$ with $x^*=x^*_1+x^*_2$ and
\begin{equation*}
\la x^*_1,\ox\ra+\la x^*_2,\ox\ra=\la x^*,\ox\ra=\sigma_{A\cap B}(x^*)=\sigma_A(x^*_1)+\sigma_B(x^*_2).
\end{equation*}
This implies that  $\la x^*_1,\ox\ra=\sigma_A(x^*_1)$ and $\la x^*_2,\ox\ra=\sigma_B(x^*_2)$. Thus we have $x^*_1\in N(\ox;A)$ and $x^*_2\in N(\ox; B)$, which show that $N(\ox;A\cap B)\subset N(\ox; A)+N(\ox; B)$. Observe that $N(\ox; A\cap B)=N(\ox;\O_1\cap\O_2)$, $N(\ox;A)=N(\ox;\O_1)$, and $N(\ox; B)=N(\ox; \O_2)$; hence we arrive at the inclusion ``$\subset$" in \eqref{ni}. The opposite inclusion therein is obvious. Similar arguments allow us to deduce normal cone intersection rule \eqref{ni} under assumptions in (ii) and (iii) from the corresponding assertions of Theorem~\ref{sigma intersection rule}.$\h$

Observe that \eqref{qc} reduces to the difference interiority condition \eqref{dqc} provided that one of the sets $\O_i$, say $\O_2$, is bounded. Note that there are other qualification conditions for the validity of \eqref{ni} in various space settings studied in connection with the so-called ``strong conic hull intersection property" (strong CHIP); see, e.g., \cite{Bauschke99,Ernst-Thera07,jey,Li-Ng-Pong07,bmn2}. Their detailed consideration is out of scope of this paper.\vspace*{-0.2in}

\section{Geometric Approach to Conjugate Calculus}
\setcounter{equation}{0}\vspace*{-0.1in}

In this section we develop a geometric approach, based on the set intersection rules obtained above, to easily derive some basic calculus results for Fenchel conjugates of extended-real-valued convex functions in LCTV and Banach spaces settings. Known proofs of such functional results are more involved and employ analytic arguments; cf.\ \cite{Ra,s,z}.\vspace*{-0.03in}

First we present a simple lemma relating the conjugate and epigraphical support functions.\vspace*{-0.15in}

\begin{Lemma}{\bf(Fenchel conjugates via supports to epigraphs).}\label{Fenchel epi} For any $f\colon X\to\oR$ we have the conjugate function representation
\begin{equation*}
f^*(x^*)=\sigma_{{\rm\small epi}\,f}(x^*,-1)\;\;\mbox{\rm whenever }\;x^*\in X^*.
\end{equation*}
\end{Lemma}\vspace*{-0.1in}
{\bf Proof.} It follows directly from the definitions that
\begin{equation*}
f^*(x^*)=\sup\big\{\la x^*,x\ra-f(x)\big|\;x\in\dom f\big\}=\sup\big\{\la x^*,x\ra-\lambda\big|\;(x,\lambda)\in\epi f\big\}=\sigma_{{\rm\small epi}\,f}(x^*,-1),
\end{equation*}
which therefore verifies the claimed relationship.$\h$

This lemma allows us to derive the conjugate sum rules and other conjugate calculus results from Theorem~\ref{sigma intersection rule}. In what follows we concentrate for simplicity on applying this theorem under the assumptions imposed in (ii) and (iii) therein.\vspace*{-0.12in}

\begin{Theorem}{\bf(conjugate sum rules).}\label{Fenchel sum rule} For convex functions $f,g\colon X\to\oR$ assume that:

{\bf(i)} either one of the functions $f$, $g$ is continuous at some point $\ox\in\dom f\cap\dom g$,

{\bf(ii)} or $X$ is Banach, $f$ and $g$ are l.s.c., and $\R^+(\dom f-\dom g)$ is a closed subspace of $X$.\\[1ex]
Then we have the conjugate to the sum $f+g$ is represented by
\begin{equation}\label{Fenchelsum}
(f+g)^*(x^*)=\big(f^*\oplus g^*\big)(x^*)\;\;\mbox{\rm for all }\;x^*\in X^*.
\end{equation}
Moreover, in both cases the infimum in $(f^*\oplus g^*)(x^*)$ is attained, i.e., for any $x^*\in\dom(f+g)^*$, there are $x^*_1,x^*_2\in X^*$ such that $x^*=x^*_1+x^*_2$ and
\begin{equation*}
(f+g)^*(x^*)=f^*(x^*_1)+g^*(x^*_2).
\end{equation*}
\end{Theorem}\vspace*{-0.12in}
{\bf Proof.} Fix $x^*\in X^*$ and easily observe that $(f+g)^*(x^*)\le(f^*\oplus g^*)(x^*)$. Indeed, for any $x^*_1, x^*_2\in X^*$ with $x^*_1+x^*_2=x^*$ we have
\begin{align*}
f^*(x^*_1)+g^*(x^*_2)&=\sup\big\{\la x^*_1,x\ra-f(x)\big|\;x\in X\big\}+\sup\big\{\la x^*_2,x\ra-g(x)\big|\;x\in X\big\}\\
&\ge\sup\big\{\la x^*_1,x\ra-f(x)+\la x^*_2,x\ra-g(x)\big|\;x\in X\big\}\\
&=\sup\big\{\la x^*,x\ra-(f+g)(x)\big|\;x\in X\big\}=(f+g)^*(x^*).
\end{align*}
Taking the infimum with respect to all such $x^*_1,x^*_2$ justifies the claimed inequality.\vspace*{-0.03in}

Let us verify the opposite inequality in \eqref{Fenchelsum} for $x^*\in\dom(f+g)^*$ in each case (i) and (ii).\\[1ex]
{\bf(i)} Define the convex subsets $\O_1,\O_2$ of $X\times\R\times\R$ by
\begin{align*}
\Omega_1:=\big\{(x,\lambda_1,\lambda_2)\in X\times\R\times\R\big|\;\lambda_1\ge f(x)\big\},\;\Omega_2:=\big\{(x,\lambda_1,\lambda_2)\in X\times\R\times\R\big|\;\lambda_2\ge g(x)\big\}
\end{align*}
and show that $(\i\Omega_1)\cap\Omega_2\ne\emp$ if $f$ is continuous at some $\ox\in\dom f\cap\dom g$. Indeed,
\begin{equation*}
\i\Omega_1=\big\{(x,\lambda_1,\lambda_2)\in X\times\mathbb{R}\times\mathbb{R}\big|\;x\in{\rm int}(\dom f),\;\lambda_1>f(x)\big\}.
\end{equation*}
It follows from the continuity of $f$ at $\ox$ that $\ox\in{\rm int}(\dom f)$. Letting $\bar{\lambda}_1:=f(\ox)+1$ and $\bar{\lambda}_2:=g(\ox)+1$ yields $(\ox,\bar{\lambda}_1,\bar{\lambda}_2)\in\i\Omega_1$ and $(\ox,\bar{\lambda}_1,\bar{\lambda}_2)\in\Omega_2$, and hence $(\i\Omega_1)\cap\Omega_2\ne\emp$. Similarly to Lemma~\ref{Fenchel epi} we have the representation
\begin{equation*}
(f+g)^*(x^*)=\sigma_{\Omega_1\cap\Omega_2}(x^*,-1,-1).
\end{equation*}
Then Theorem~\ref{sigma intersection rule}(ii) gives us triples $(x_1^*,-\alpha_1,-\alpha_2),(x_2^*-\beta_1,-\beta_2)\in X^*\times\R\times\R$
such that $(x^*,-1,-1)=(x_1^*,-\alpha_1,-\alpha_2)+(x_2^*,-\beta_1,-\beta_2)$ and
\begin{equation*}
(f+g)^*(x^*)=\sigma_{\Omega_1\cap\Omega_2}(x^*,-1,-1)=\sigma_{\Omega_1}(x^*_1,-\alpha_1,-\alpha_2)+\sigma_{\Omega_2}(x^*_2,-\beta_1,-\beta_2).
\end{equation*}
If $\alpha_2\ne 0$, then $\sigma_{\Omega_1}(x_1^*,-\alpha_1,-\alpha_2)=\infty$, which is not possible since $x^*\in\dom(f+g)^*$, and so $\alpha_2=0$. Similarly we get $\beta_1=0$. Employing Lemma~\ref{Fenchel epi} again and taking into account the structures of $\Omega_1$ and $\Omega_2$ tell us that
\begin{align*}
(f+g)^*(x^*)&=\sigma_{\Omega_1\cap\Omega_2}(x^*,-1,-1)=\sigma_{\Omega_1}(x^*_1,-1,0)+\sigma_{\Omega_2}(x^*_2,0,-1)\\
&=\sigma_{{\rm\small epi}\,f}(x^*,-1)+\sigma_{{\rm\small epi}\,g}(x^*_2,-1)=f^*(x^*_1)+g^*(x^*_2)\ge(f^*\oplus g^*)(x^*),
\end{align*}
which verifies \eqref{Fenchelsum} and also the last statement of the theorem.\\[1ex]
{\bf(ii)} To prove the theorem under the conditions in (ii), observe first that the sets $\O_1,\O_2$ above are closed by the l.s.c.\ assumption on $f,g$ and then  check that
\begin{equation}\label{ab-dom}
\R^+(\Omega_1-\Omega_2)=\R^+\big(\dom f-\dom g\big)\times\R\times\R.
\end{equation}
Indeed, consider $u\in\mathbb{R}^+(\Omega_1-\Omega_2)$ and find $t>0$, $v\in(\Omega_1-\Omega_2)$ with $u=tv$; therefore $v=(x_1,\lambda_1,\lambda_2)-(x_2,\gamma_1,\gamma_2)$ with $(x_1,\lambda_1,\lambda_2)\in\Omega_1$ and $(x_2,\gamma_1,\gamma_2)\in\Omega_2$. Note that $x_1\in\dom f$ and $x_2\in\dom g$ due to $f(x_1)\le\lambda_1<\infty$ and $g(x_2)\le\gamma_2<\infty$. Hence we get
\begin{equation*}
tv=t(x_1-x_2,\lambda_1-\gamma_1,\lambda_2-\gamma_2)\in\mathbb{R}^+(\dom f-\dom g)\times\mathbb{R}\times\mathbb{R}.
\end{equation*}\vspace*{-0.25in}

To verify the opposite inclusion in \eqref{ab-dom}, fix $x\in\mathbb{R}^+(\dom f-\dom g)\times\mathbb{R}\times\mathbb{R}$ and find $t>0$, $x_1\in\dom f$, $x_2\in \dom g$, $\gamma_1,\gamma_2\in\mathbb{R}$, and $\lm_1,\lm_2$ with
\begin{align*}
x&=(t(x_1-{x_2}),\gamma_1,\gamma_2)=t(x_1-x_2,\lambda_1,\lambda_2)\\
&= t\big[\big(x_1,f(x_1),\lambda_2+g(x_2)\big)-(x_2,-\lambda_1+f(x_1),g(x_2)\big)\big],
\end{align*}
which readily yields $x\in t(\Omega_1-\Omega_2)\subset\mathbb{R}^+(\Omega_1-\Omega_2)$. Applying now Theorem~\ref{sigma intersection rule}(iii), we arrive at both conclusions of the theorem under (ii) and thus complete the proof. $\h$

The next result derive geometrically the chain rules for Fenchel conjugates.\vspace*{-0.12in}

\begin{Theorem}{\bf(conjugate chain rules).}\label{Fenchel chain rule} Let $A\colon X\rightarrow Y$ be a linear continuous mapping, and let $g\colon Y\rightarrow \mathbb{\overline{R}}$ be a convex function. Assume that:

{\bf(i)} either $g$ is finite and continuous at some point of $AX$,

{\bf(ii)} or $X$ and $Y$ are Banach, $g$ is l.s.c., and $\R^+(\dom g-AX)$ is a closed subspace of $Y$.\\[1ex]
Then the conjugate of the composition $g\circ A$ is represented by
\begin{equation}\label{con-chain}
(g\circ A)^*(x^*)=\inf\big\{g^*(y^*)\big|\;y^*\in(A^*)^{-1}(x^*)\big\}.
\end{equation}
Moreover, in both cases we have that for any $x^*\in\dom(g\circ A)^*$ there is $y^*\in(A^*)^{-1}(x^*)$ with
\begin{equation*}
(g\circ A)^*(x^*)=g^*(y^*).
\end{equation*}
\end{Theorem}\vspace*{-0.1in}
{\bf Proof}. Fix any $y^*\in(A^*)^{-1}(x^*)$ and deduce the following relationships from the definitions:
\begin{align*}
g^*(y^*)&=\sup\big\{\la y^*,y\ra-g(y)\big|\;y\in Y\big\}\\
&\ge\sup\big\{\la y^*,A(x)\ra-g(A(x))\big|\;x\in X\big\}\\
&=\sup\big\{\la A^*y^*,x\ra-(g\circ A)(x)\big|\;x\in X\big\}\\
&=\sup\big\{\la x^*,x\ra-(g\circ A)(x)\big|\;x\in X\big\}=(g\circ A)^*(x^*).
\end{align*}
This readily implies the inequality ``$\le$" in \eqref{con-chain}. To verify the opposite inequality therein, take $x^*\in\dom(g\circ A^*)$ and construct the following convex sets:
\begin{align*}
\Omega_1:=\gph A\times\R\subset X\times Y\times\R\;\;\mbox{\rm and }\;\Omega_2:=X\times\epi g\subset X\times Y\times\R.
\end{align*}
It follows directly from the above constructions that
\begin{equation*}
(g\circ A)^*(x^*)=\sigma_{\Omega_1\cap\Omega_2}(x^*,0,-1)<\infty.
\end{equation*}
Let us further proceed with proof in each case (i) and (ii).\\[1ex]
{\bf (i)} Under the assumption in (i) we clearly have
\begin{equation*}
{\rm int}\,\Omega_2=\big\{(x,y,\lambda)\in X\times Y\times\R\big|\;x\in X,\;y\in{\rm int}(\dom g),\;\lambda>g(y)\big\}
\end{equation*}
and easily observe that $\Omega_1\cap(\i\Omega_2)\ne\emp$ for the sets $\O_1,\O_2$ above. Then Theorem~\ref{sigma intersection rule}(ii) tells us that there are triples $(x^*_1,y^*_1,\alpha_1),(x^*_2, y^*_2,\alpha_2)\in X^*\times Y^*\times\R$ satisfying
\begin{equation*}
(x^*,0,-1)=(x^*_1,y^*_1,\alpha_1)+(x^*_2,y^*_2,\alpha_2)\;\mbox{\rm and }\sigma_{\Omega_1\cap\Omega_2}(x^*,0,-1)=\sigma_{\Omega_1}(x^*_1,y^*_1,\alpha_1)+\sigma_{\Omega_2}(x^*_2,y^*_2,\alpha_2).
\end{equation*}
It follows from the structures of $\Omega_1,\Omega_2$ that $\alpha_1=0$ and $x^*_2=0$. This yields
\begin{equation*}
\sigma_{\Omega_1\cap\Omega_2}(x^*,0,-1)=\sigma_{\Omega_1}(x^*,-y^*_2,0)+\sigma_{\Omega_2}(0,y^*_2,-1)
\end{equation*}
for some $y^*_2\in Y^*$. Thus we arrive at the representations
\begin{align*}
\sigma_{\Omega_1\cap\Omega_2}(x^*,0,-1)&=\sup\big\{\la x^*,x\ra-\la y^*_2,A(x)\ra\big|\;x\in X\ra\big\}+\sigma_{{\rm\small epi}\,g}(y^*_2,-1)\\
&=\sup\big\{\la x^*-A^*y_2^*,x\ra\big |\; x\in X\big\}+g^*(y^*_2),
\end{align*}
which allow us to conclude that $x^*=A^*y_2^*$ and therefore
\begin{equation*}
\sigma_{\Omega_1\cap\Omega_2}(x^*,0,-1)=g^*(y^*_2)\ge\inf\big\{g^*(y^*)\big|\;y^*\in(A^*)^{-1}(x^*)\big\}.
\end{equation*}
This justifies both statements of the theorem in case (i).\\[1ex]
{\bf (ii)} It is easy to check the equality
\begin{equation*}
\R^{+}(\Omega_1-\Omega_2)=X\times\R^{+}(AX-\dom g)\times\R.
\end{equation*}
Then we apply Theorem~\ref{sigma intersection rule}(iii) and get both claimed statement in this case as well.$\h$

Finally in this section, we present a simple geometric proof of the Fenchel conjugate representation for the {\em maximum} of two convex functions
$f,g\colon X\to\oR$ defined by
\begin{equation}\label{max}
(f\vee g)(x):=\max\big\{f(x),g(x)\big\},\quad x\in X.
\end{equation}

In the next theorem we use the convention that $0f:=\delta_{\dom f}$ and similarly for $g$; see \cite{Ra1}.\vspace*{-0.12in}

\begin{Theorem}{\bf(calculating conjugates of maximum functions).}\label{Fenchel max rule} Given convex functions $f,g\colon X\to\oR$, suppose that either the assumptions in {\bf(i)} or the assumptions in {\bf(ii)} of Theorem~{\rm\ref{Fenchel sum rule}} are satisfied. Then we have the representation
\begin{equation}\label{Fenchelmax}
(f\vee g)^*(x^*)=\inf_{\lambda\in[0,1]}\big[\lambda f+(1-\lambda)g\big]^*(x^*).
\end{equation}
If furthermore $(f\vee g)^*(x^*)\in\R$, then the minimum in \eqref{Fenchelmax} is achieved.
\end{Theorem}\vspace*{-0.1in}
{\bf Proof.} Let us first check that the inequality ``$\le$" always holds in \eqref{Fenchelmax} whenever $\lm\in[0,1]$. Indeed, it follows directly from the definitions that
\begin{align*}
\big[\lambda f+(1-\lambda)g\big]^*(x^*)&=\sup_{x\in X}\big[\la x^*, x\ra-\lambda f(x)-(1-\lambda)g(x)\big]\\
&\ge\sup_{x\in X}\big[\la x^*,x\ra-\lambda(f\vee g)(x)-(1-\lambda)(f\vee g)(x)\big]\\
&=\sup_{x\in X}\big[\la x^*,x\ra-(f\vee g)(x)\big]=(f\vee g)^*(x^*),\quad x^*\in X^*.
\end{align*}
To verify the opposite inequality, observe that $\epi(f\vee g)=\epi f\cap\epi g$, and hence we deduce from Lemma~\ref{Fenchel epi} the relationships
\begin{equation*}
(f\vee g)^*(x^*)=\sigma_{{\rm epi}(f\vee g)}(x^*,-1)=\sigma_{\Omega_1\cap\Omega_2}(x^*,-1)\;\;\mbox{\rm with }\;\Omega_1:=\epi f,\;\Omega_2=\epi g.
\end{equation*}
It follows from Theorem~\ref{sigma intersection rule} under the assumptions in either (i) or (ii) that
\begin{equation*}
(f\vee g)^*(x^*)=\big[\sigma_{\Omega_1}\oplus\sigma_{\Omega_2}\big](x^*,-1),\quad x^*\in X^*.
\end{equation*}
Observe that the continuity of $f$ at $\ox\in\dom f\cap\dom g$ yields $({\rm int}\,\O_1)\cap\O_2\ne\emp$. In the case where $\R^+(\dom f-\dom g)$ is a closed subspace we have that the set $\R^+(\O_1-\O_2)=\R^+(\dom f-\dom g)\times\R$ is a closed subspace as well.\vspace*{-0.03in}

Suppose now that $(f\vee g)^*(x^*)\in\R$, i.e., $x^*\in\dom(f\vee g)^*$. Then Theorem~\ref{sigma intersection rule} gives us pairs $(x^*_1,-\lambda_1),(x^*_2,-\lm_2)\in X^*\times\R$ such that $x^*=x^*_1+x^*_2$, $1=\lambda_1+\lambda_2$, and
\begin{equation*}
(f\vee g)^*(x^*)=\sigma_{{\rm\small epi}(f\vee g)}(x^*,-1)=\sigma_{\Omega_1\cap\Omega_2}(x^*,-1)=\sigma_{\Omega_1}(x^*_1,-\lambda_1)+\sigma_{\Omega_2}(x^*_2,-\lambda_2).
\end{equation*}
Note that if either $\lambda_1<0$ or $\lambda_2<0$, then $\sigma_{\Omega_1}(x^*_1,-\lambda_1)=\infty$ or $\sigma_{\Omega_2}(x^*_1,-\lambda_1)=\infty$, which is a contradiction. Thus $\lambda_1,\lm_2\ge 0$. In the case where $\lambda_1,\lm_2>0$ it follows that
\begin{align*}
(f\vee g)^*(x^*)&=\lambda_1\sigma_{\Omega_1}\Big(\frac{x^*_1}{\lambda_1},-1\Big)+\lambda_2\sigma_{\Omega_2}\Big(\frac{x^*_2}{\lambda_2}, -1\Big)=\lambda_1f^*\Big(\frac{x^*_1}{\lambda_1}\Big)+\lambda_2g^*\Big(\frac{x^*_2}{\lambda_2}\Big)\\
&=(\lambda_1f)^*(x^*_1)+(\lambda_2g)^*(x^*_2)\;\;\mbox{\rm with }\;x^*\in\dom(f\vee g)^*.
\end{align*}
Furthermore, we obviously have the estimate
\begin{equation*}
(\lambda_1f)^*(x^*_1)+(\lambda_2g)^*(x^*_2)\ge\big[\lambda_1 f+(1-\lambda_1)g\big]^*(x^*)\ge\inf_{\lambda\in[0,1]}\big[\lambda f+(1-\lambda)g\big]^*(x^*).
\end{equation*}
Plugging there $\lambda_1=0$ gives us $\lambda_2=1$, and hence
\begin{align*}
\sigma_{\Omega_1}(x^*_1,-\lambda_1)+\sigma_{\Omega_2}(x^*_2,-\lambda_2)&=\sigma_{\dom f}(x^*_1)+g^*(x^*_2)\\
&=\delta^*_{\dom f}(x^*_1)+g^*(x^*_2)\ge\big[\delta_{\dom f}+g\big]^*(x^*)\\
&\ge\inf_{\lambda\in[0,1]}\big[\lambda f+(1-\lambda)g\big]^*(x^*).
\end{align*}
Since the latter inequality also holds if $(f\vee g)^*(x^*)=\infty$, we complete the proof. $\h$\vspace*{-0.2in}

\section{Fenchel Conjugate of Optimal Value Functions}
\setcounter{equation}{0}\vspace*{-0.1in}

In this section we start our study of a class of extended-real-valued functions, which plays a highly important role in many aspects of variational analysis, optimization, and their applications; see, e.g., \cite{m-book1,RockWets-VA,thib} and the references therein. Such objects, known as {\em marginal} or {\em optimal value functions} are constructed in the following way. Given a set-valued mapping $F\colon X\tto Y$ and an extended-real-valued function $\ph\colon X\times Y\to\oR$, define
\begin{equation}\label{marg}
\mu(x):=\inf\{\ph(x,y)\big|\;y\in F(x)\big\},\quad x\in X.
\end{equation}
We always assume that $\mu(x)>-\infty$ for all $x\in X$ and observe that the defined function \eqref{marg} is convex provided that $\ph$ is convex and $F$ is convex-graph; see, e.g., \cite[Proposition~1.50]{bmn}.\vspace*{-0.03in}

The main goal of this short section is to derive a new formula for calculating the Fenchel conjugate of the optimal value function \eqref{marg} expressed in terms of the initial data $\ph$ and $F$. We continue with subdifferential calculus rules for \eqref{marg} in Section~8.\vspace*{-0.12in}

\begin{Theorem}{\bf(calculating Fenchel conjugates of optimal value functions).}\label{prop:opt_fenchel} Let $\ph\colon X\times Y\to\oR$ be a convex function, and let the graph of the mapping $F\colon X\tto Y$ be a convex subset of $X\times Y$. Then we have
\begin{equation}\label{marg1}
\mu^*(x^*)=\big(\ph+\delta_{{\rm\small gph}\,F}\big)^*(x^*,0)\;\;\mbox{\rm whenever }\;x^*\in X^*.
\end{equation}
Furthermore, the following refined representation
\begin{equation}\label{marg2}
\mu^*(x^*)=\big(\ph^*\oplus\sigma_{{\rm\small gph}\,F}\big)(x^*,0),\quad x^*\in X^*,
\end{equation}
holds provided the validity of one of the assumption groups below:

{\bf(i)} either $\ph$ is finite and continuous at some point $(\ox,\oy)\in\gph F$,

{\bf(ii)} or $X$ and $Y$ are Banach, $\ph$ is l.s.c., $F$ is of closed graph, and the set $\R^{+}(\dom\ph-\gph F)$ is a closed subspace of $X\times Y$.
\end{Theorem}\vspace*{-0.1in}
{\bf Proof.} Fix $x^*\in X^*$ and $x\in\dom\mu$. It follows from definition \eqref{marg} that
\begin{equation*}
\la x^*,x\ra-\mu(x)\ge\la x^*,x\ra-\ph(x,y)\;\text{ whenever }\;y\in F(x).
\end{equation*}
This clearly implies that for all $(x,y)\in X\times Y$ we have
\begin{align*}
\mu^*(x^*)=&\sup\big\{\la x^*,x\ra-\mu(x)\big|\;x\in\dom\mu\big\}\ge\la x^*,x\ra-\mu(x)\\
\ge&\la(x^*,0),(x,y)\ra-\big(\ph+\delta_{{\rm\small gph}\,F}\big)(x,y).
\end{align*}
Getting these things together gives us the relationships
\begin{align*}
\mu^*(x^*)\ge\sup\big\{\la(x^*,0),\,(x,y)\ra-\big(\ph+\delta_{{\rm\small gph}\,F}\big)(x,y)\big|\;(x,y)\in X\times Y\big\}=\big(\ph+\delta_{{\rm\small gph}\, F}\big)^*(x^*,0),
\end{align*}
which justify the inequality ``$\ge$" in \eqref{marg1}. To verify the opposite inequality therein, fix $\ve>0$ and for any $x\in\dom\mu$ find $\oy\in F(x)$ such that $\ph(x,\oy)<\mu(x)+\epsilon$. Then we have
\begin{align*}
\la x^*,x\ra-\big(\mu(x)+\epsilon\big)<&\la x^*,x\ra-\ph(x,\oy)\le\sup\big\{\la x^*,x\ra-\ph(x,y)\big|\;y\in F(x)\big\}\\
\le&\sup\big\{\la(x^*,0),\,(x,y)\ra-\big(\ph+\delta_{{\rm\small gph}\,F}\big)(x,y)\big|\;(x,y)\in X\times Y\big\}\\
=&\big(\ph+\delta_{{\rm\small gph}\,F}\big)^*(x^*,0)\;\text{ whenever }\;x^*\in X^*.
\end{align*}
Since these relationships hold for all $x\in\dom\mu$ and all $\epsilon>0$, we conclude that
\begin{equation*}
\mu^*(x^*)\le\big(\ph+\delta_{{\rm\small gph}\,F}\big)^*(x^*,0),\quad x^*\in X^*,
\end{equation*}
which therefore justifies the first representation \eqref{marg1} of the theorem.\vspace*{-0.03in}

To derive the second representation \eqref{marg2} for the conjugate of the optimal value function, it remains applying the conjugate sum rule from Theorem~\ref{Fenchel sum rule} to the sum in \eqref{marg1} with taking into account that $\dd^*_{{\rm\small gph}\,F}(x^*,0)=\sigma_{{\rm\small gph}\,F}(x^*,0)$ for all $x^*\in X^*$. $\h$\vspace*{-0.15in}

\section{Geometric Approach to Convex Subdifferential Calculus}
\setcounter{equation}{0}\vspace*{-0.1in}

This section presents a geometric design of major rules of convex subdifferential calculus (including some new results) that is based on applications of the normal cone intersection rule taken from Theorem~\ref{nir}. Having in mind a clear illustration of the developed geometric approach, we mainly confine ourselves to the applications of the easily formulated qualification conditions in (ii) and (iii) therein. Let us start with the fundamental sum rule in both LCTV and Banach space settings.\vspace*{-0.15in}

\begin{Theorem}{\bf(subdifferential sum rule).}\label{sub sum rule} Given convex functions $f,g\colon X\to\oR$, we have the subdifferential sum rule
\begin{equation}\label{sumrule}
\partial(f+g)(\bar{x})=\partial f(\ox)+\partial g(\ox)\;\;\mbox{\rm for all }\;\ox\in\dom f\cap\dom g
\end{equation}
provided the validity of the following:\\
{\bf(i)} either $f$ is finite and continuous at some point in $\dom f\cap\dom g$,\\
{\bf(ii)} or $X$ is a Banach space, $f$ and $g$ are l.s.c., and $\R^{+}(\dom f-\dom g)$ is a subspace of $X$.
\end{Theorem}\vspace*{-0.12in}
{\bf Proof.} Let $\ox\in\dom f\cap\dom g$ be fixed for the entire proof. Since the inclusion ``$\supset$" in \eqref{sumrule} can be easily checked by the definition, we now concentrate on proving the opposite inclusion therein. Pick any $x^*\in\partial(f+g)(\bar{x})$ and show how the geometric results of Theorem~\ref{nir} can be used in verifying $x^*\in\partial f(\ox)+\partial g(\ox)$. Having
\begin{equation*}
\langle x^*,x-\bar{x}\rangle\le(f+g)(x)-(f+g)(\bar{x})\;\mbox{\rm for all }\;x\in X,
\end{equation*}
define the following convex subsets of $X\times\R\times\R$:
\begin{align*}
\Omega_1:= &\big\{(x,\lambda_1,\lambda_2)\in X\times\mathbb{R}\times\mathbb{R}\big|\;\lambda_1\ge f(x)\big\},\\
\Omega_2:= & \big\{(x,\lambda_1,\lambda_2)\in X\times\mathbb{R}\times\mathbb{R}\big|\;\lambda_2\ge g(x)\big\}.
\end{align*}
It follows from the definition that $(x^*,-1,-1)\in N((\bar{x},f(\bar{x}),g(\bar{x}));\Omega_1\cap\Omega_2)$. The application of Theorem~\ref{nir} under an appropriate qualification condition yields
\begin{equation}\label{eq:sumint}
(x^*,-1,-1)\in N\big((\bar{x},f(\bar{x}),g(\bar{x}));\Omega_1\big)+N\big((\bar{x},f(\bar{x}),g(\bar{x}));\Omega_2\big),
\end{equation}
which tells us therefore that
\begin{equation*}
(x^*,-1,-1)=(x^*_1,-\lambda_1,-\lambda_2)+(x^*_2,-\gamma_1,-\gamma_2)
\end{equation*}
with $(x^*_1,-\lambda_1,-\lambda_2)\in N((\bar{x},f(\bar{x}),g(\bar{x}));\Omega_1)$ and $(x^*_2,-\gamma_1,-\gamma_2)\in N((\bar{x},f(\bar{x}),g(\bar{x}));\Omega_2)$. By the construction of $\Omega_1$ and $\Omega_2$ we have $\lambda_2=\gamma_1=0$ and hence find dual elements $(x^*_1,-1)\in N((\bar{x},f(\bar{x}));\epi f) $ and $(x^*_2,-1)\in N((\bar{x},g(\bar{x}));\epi g)$ satisfying the relationships
\begin{equation*}
x^*_1\in\partial f(\bar{x}),\quad x^*_2\in\partial g(\bar{x}),\;\;\mbox{\rm and }\;x^*=x^*_1+x^*_2.
\end{equation*}
This shows that $x^*\in\partial f(\bar{x})+\partial g(\bar{x})$, and thus \eqref{sumrule} holds, provided that the corresponding conditions of Theorem~\ref{nir} are satisfied under the assumptions imposed in the theorem.\vspace*{-0.03in}

To this end, we easily observe that (i) yields $(\i\Omega_1)\cap\Omega_2\ne\emp$ and that (ii) ensures the closedness of the subspace
\begin{eqnarray*}
\R^+(\Omega_1-\Omega_2)=\R^+(\dom f-\dom g)\times\R\times\R
\end{eqnarray*}
for the sets $\O_1,\O_2$ defined above. $\h$

Next we employ the geometric approach to obtain a chain rule for convex subgradients under different qualification conditions.\vspace*{-0.12in}

\begin{Theorem}{\bf(subdifferential chain rule).}\label{sub chainrule}
Let $A\colon X\rightarrow Y$ be a continuous linear mapping, and let $g\colon Y\rightarrow \mathbb{\overline{R}}$ be convex. Suppose that: \\[1ex]
{\bf(i)} either $g$ is finite and continuous at some point of $AX$.\\
{\bf(ii)} or $X$ is Banach, $g$ is l.s.c., and $\R^+(AX-\dom g)$ is a closed subspace of $X$.\\[1ex]
Then we have the subdifferential chain rule
\begin{equation}\label{eq3.0.6}
\partial(g\circ A)(\bar{x})=A^*\partial g(A\ox):=\big\{A^*y^*\big|\;y^*\in\partial g(A\ox)\big\}\;\;\mbox{\rm for all }\;\ox\in{\rm dom}(g\circ A).
\end{equation}
\end{Theorem}\vspace*{-0.2in}
{\bf Proof.} First we verify the inclusion ``$\subset$" in \eqref{eq3.0.6} under the validity of (i) and (ii). Pick any $y^*\in\partial(g\circ A)(\bar{x})$ and form the subsets of $X\times Y\times\mathbb{R}$ by
\begin{align*}
\Omega_1:=(\gph A)\times\mathbb{R}\;\;\mbox{\rm and }\;\Omega_2:= X\times(\epi g).
\end{align*}
If (i) holds, then it follows from the proof of Theorem~\ref{Fenchel chain rule}(i) that $\Omega_1\cap(\i\Omega_2)\ne\emp$. Using the definitions of subgradients and normals of convex analysis, we easily conclude that $(y^*,0,-1)\in N((\bar{x},\bar{y},\bar{z});\Omega_1\cap\Omega_2)$ with $\bar{z}:=g(\bar{y})$ and $\oy:=A\ox$. The intersection rule from Theorem~\ref{nir}(ii) tells us that
\begin{equation*}
(y^*,0,-1)\in N\big((\bar{x},\bar{y},\bar{z});\Omega_1\big)+ N\big((\bar{x},\bar{y},\bar{z});\Omega_2\big),
\end{equation*}
which ensures the existence of $(y^*,-w^*)\in N((\bar{x},\bar{y});\gph A)$ and $(w^*,-1)\in N((\bar{y},\bar{z});\epi g)$ satisfying $(y^*,0,-1)=(y^*,-w^*,0)+(0,w^*,-1)$. We can directly check that
\begin{equation*}
N\big((\bar{x},\bar{y});\gph A\big)=\big\{(x^*,y^*)\in X^*\times Y^*\big|\;x^*=-A^*y^*\big\}
\end{equation*}
and therefore deduce from the above that
\begin{equation*}
y^*= A^* w^*\;\;\mbox{\rm and }\;w^*\in\partial g(\bar{y}).
\end{equation*}
This gives us $y^*\in A^*\partial g(\bar{y})$ and hence justifies the claimed inclusion ``$\subset$" in \eqref{eq3.0.6}.\vspace*{-0.03in}

To verify the opposite inclusion in \eqref{eq3.0.6}, fix $y^*\in\partial g(A(\bar{x}))$ and deduce from the definitions of convex subgradients and adjoint linear operators that
\begin{align*}
\langle y^*,Ax-A\bar{x}\rangle\le g(Ax)-g(A\bar{x})\Longleftrightarrow\langle A^*y^*,x-\bar{x}\rangle\le(g\circ A)(x)-(g\circ A)(\bar{x})
\end{align*}
for every $x\in X$. This implies that $ A^*\partial g(\bar{y})\subset\partial(g\circ A)(\bar{x})$, which completes the proof of the theorem under the assumptions in (i)\vspace*{-0.03in}

If now the assumptions in (ii) are satisfied, then we can easily see that
\begin{equation*}
\R^{+}(\Omega_1-\Omega_2)=X\times\R^{+}(AX-\dom g)\times\R.
\end{equation*}
It remains to employ Theorem~\ref{nir}(iii) in this setting by using the arguments above. $\h$

The final result of this section employs the geometric approach to derive formulas for subdifferentiation of maximum of convex functions. For simplicity we consider the case of two functions; the maximum of finitely many ones can be deduced by induction.\vspace*{-0.12in}

\begin{Theorem}{\bf(subdifferentiation of maximum functions).}\label{sub-max} Let $\bar{x}\in \dom f\cap \dom g$, where the functions $f, g\colon X\to\oR$ are convex. Then the following assertions hold:\\[1ex]
{\bf(a)} If $f(\ox)>g(\ox)$ and $g$ is upper semicontinuous at $\ox$, then
\begin{equation*}
\partial(f\vee g)(\ox)=\partial f(\ox).
\end{equation*}
{\bf(b)} If $f(\ox)<g(\ox)$ and $f$ is upper semicontinuous at $\ox$, then
\begin{equation*}
\partial(f\vee g)(\ox)=\partial g(\ox).
\end{equation*}
{\bf(c)} If $f(\ox)=g(\ox)$, then
\begin{equation}\label{max-rule}
\partial (f\vee g)(\ox)=\co[\partial f(\ox)\cup \partial g(\ox)].
\end{equation}
provided that either the assumptions in {\bf(i)} or the assumptions in {\bf(ii)} of Theorem~{\rm\ref{Fenchel sum rule}} are satisfied. Here the symbol `${\rm co}$' stands for the convex hull of a set.
\end{Theorem}\vspace*{-0.1in}
{\bf Proof.} We obviously have ${\rm epi}(f\vee g)=\epi f\cap\epi g$, and thus it follows from the subgradient definition that $x^*\in\partial (f\vee g)(\ox)$ if and only if
\begin{equation}\label{nor-max}
(x^*,-1)\in N\big((\ox,(f\vee g)(\ox));\epi (f\vee g)\big)=N\big((\ox,(f\vee g)(\ox));\epi f\cap\epi g\big).
\end{equation}
Denote $\bar{\lambda}:=(f\vee g)(\ox)$. Under the assumption given in (a), it is easy to check that $(\ox,\bar\lambda)\in{\rm int}(\epi g)$  due to the imposed upper semicontinuity assumption on $g$. This allows us to exclude the latter function from the intersection in \eqref{nor-max} and arrive at the conclusion. The proof of assertion (b) is the same.\vspace*{-0.03in}

It remains to verify formula \eqref{max-rule} under the assumptions given in (c). To furnish this, we only need to apply Theorem~\ref{nir} and proceed similarly to the proof of Theorem~\ref{Fenchel max rule}. $\h$\vspace*{-0.2in}

\section{Subgradients of Optimal Value Functions}
\setcounter{equation}{0}\vspace*{-0.1in}

In this section we return to the class of (convex) optimal value/marginal functions defined in \eqref{marg}. Observing that functions \eqref{marg} are {\em intrinsically nonsmooth}, we concentrate here on calculating their subgradient sets. Note that the new results obtained below are significantly different from those known earlier for nonconvex marginal functions (see, e.g., \cite[Theorems~1.108 and 3.38]{m-book1} with the references therein) and their specifications to the convex case. We present an application of the obtained results for marginal functions to deriving another {\em chain rule} for mappings defined on LCTV spaces with ordered values. \vspace*{-0.03in}

Given the optimal value function $\mu(x)$ which composed in \eqref{marg}, define the corresponding {\em solution map} $M\colon X\tto Y$ associated with \eqref{marg} by
\begin{equation}\label{sol}
M(x):=\big\{y\in F(x)\big|\;\mu(x)=\ph(x,y)\big\},\quad x\in X.
\end{equation}
The first result of this section allows us to represent the subdifferential of the optimal value function via the one for the sum of two convex functions taken from definition \eqref{marg}. \vspace*{-0.12in}

\begin{Proposition}{\bf(sum representation for subgradients of optimal value functions).}\label{marginal sum} Considering \eqref{marg}, let $\ph\colon X\times Y\to\oR$ be convex, and let $F\colon X\tto Y$ be convex-graph. Then for any $\ox\in X$ and $\oy\in M(\ox)$ we have the subdifferential representation
\begin{align*}
\partial\mu(\ox)=\big\{x^*\in X^*\big|\;(x^*,0)\in\partial\big(\ph+\delta_{{\rm\small gph}\,F}\big)(\ox,\oy)\big\}.
\end{align*}
\end{Proposition}\vspace*{-0.1in}
{\bf Proof.} Fix a subgradient $x^*\in\partial\mu(\ox)$ and get by \eqref{sub} that
\begin{equation}\label{subestimate1}
\la x^*,x-\ox\ra\le\mu(x)-\mu(\ox)\;\;\mbox{\rm for all }\;x\in X.
\end{equation}
It follows from the definition in \eqref{marg} that
\begin{align*}
\la x^*,x-\ox\ra\le\mu(x)-\ph(\ox,\oy)\le\ph(x,y)-\ph(\ox,\oy)\;\;\mbox{\rm whenever }\;y\in F(x).
\end{align*}
Using further the indicator function associated with $\gph F$ tells us that
\begin{equation}\label{subestimate}
\la x^*, x-\ox\ra\le\mu(x)-\ph(\ox,\oy)\le\big(\ph+\delta_{{\rm\small gph}\,F}\big)(x,y)-\big(\ph+\delta_{{\rm\small gph}\,F}\big)(\ox,\oy)
\end{equation}
for all $(x,y)\in X\times Y$, which yields the inclusion $(x^*,0)\in\partial(\ph+\delta_{{\rm\small gph}\,F})(\ox,\oy)$.\vspace*{-0.03in}

To verify the opposite inclusion, pick $(x^*,0)\in\partial(\ph+\delta_{{\rm\small gph}\,F})(\ox,\oy)$. Then \eqref{subestimate} is satisfied for all $(x,y)\in X\times Y$. Therefore
\begin{align*}
\la x^*,x-\ox\ra\le\ph(x,y)-\ph(\ox,\oy)=\ph(x,y)-\mu(\ox)\;\;\mbox{\rm whenever }\;y\in F(x).
\end{align*}
Taking the infimum on the right-hand side above with respect to $y\in F(x)$ gives us \eqref{subestimate1}, and so we arrive at $x^*\in\partial\mu(\ox)$ and complete the proof of the proposition. $\h$

The next theorem provides a precise calculation of the subdifferential of the optimal value functions entirely in terms of the generalized differential constructions for the initial data $\ph$ and $F$ of \eqref{marg}, namely via the subdifferential of $\ph$ and the coderivative of $F$. Some particular case of the obtained representation can be found in \cite{AnYen,bmn,bmn1}.\vspace*{-0.12in}

\begin{Theorem}{\bf(precise subdifferential calculation for optimal value functions).}\label{AnYenThm42} Let $\ph\colon X\times Y\to\oR$ be a convex function, and let $F\colon X\tto Y$ be a set-valued mapping with convex graph. Assume the following:\\[1ex]
{\bf (i)} either $\ph$ is finite and continuous at some point in $\gph F$,\\
{\bf (ii)} or $X$ is Banach, $\ph$ is l.s.c., the graph of $F$ is closed, and $\R^{+}(\dom\ph-\gph F)$ is a closed subspace of $X\times Y$.\\[1ex]
Then for any $\ox\in X$ and $\oy\in M(\ox)$ we have the subdifferential representation
\begin{equation}\label{marg-sub}
\partial\mu(\ox)=\bigcup_{(x^*,y^*)\in\partial\ph(\ox,\oy)}\big[x^*+D^*F(\ox,\oy)(y^*)\big].
\end{equation}
\end{Theorem}\vspace*{-0.1in}
{\bf Proof.} Having in hands the subdifferential representation for the optimal value function obtained in Proposition~\ref{marginal sum}, we apply therein the subdifferential sum rule from Theorem~\ref{sub sum rule} under the corresponding conditions in (i) and (ii). The coderivative of $F$ appears there due to its definition and the fact that $\partial\dd_\O(\oz)=N(\oz;\O)$ for any convex set $\O$. $\h$\vspace*{-0.12in}

\begin{Corollary}{\bf(subdifferential of optimal value functions with parameter independent costs).}\label{remark:FundThm} Suppose that the cost function $\ph$ in \eqref{marg} does not depend on $x$. Then we have
\begin{align*}
\partial\mu(\ox)=\bigcup_{y^*\in\partial\ph(\oy)}D^*F(\ox,\oy)(y^*)\;\;\mbox{\rm for any }\;(\ox,\oy)\in\gph M.
\end{align*}
with $M$ from \eqref{sol} under the validity of one of the qualification conditions {\rm(i)} and {\rm(ii)}.
\end{Corollary}\vspace*{-0.12in}
{\bf Proof.} This follows directly from \eqref{marg-sub} with $\ph(x,y)=\ph(x)$. Indeed, in this case we obviously have $\partial\ph(\ox,\oy)=\partial\ph(\ox)$ and $x^*=0$. $\h$\vspace*{-0.12in}

\begin{Remark}{\bf(chain rules from subdifferentiation of marginal functions).}\label{cr-marg}{\rm When the mapping $F$ in \eqref{marg} is {\em single-valued}, the marginal function therein reduces to the (generalized) {\em composition} $\ph(x,F(x))$, which is the standard composition $\ph\circ F$ if $\ph=\ph(y)$. In this way the results of Theorem~\ref{AnYenThm42} and Corollary~\ref{remark:FundThm} give us some (generalized) {\em chain rules}. Consider, e.g., the case where $F(x):=Ax+b$, $A:X\to Y$ is a continuous linear operator, and $b\in Y$. Then the coderivative of $F$ at $(\ox,\oy)$ with $\oy:=A\ox+b$ is calculated by
\begin{align*}
D^*F(\ox,\oy)(y^*)=\{A^*y^*\}\;\;\mbox{\rm for all }\;y^*\in Y^*.
\end{align*}
The composition $g\circ F$ of this mapping $F$ with a convex function $g\colon Y\to\oR$ is a particular case of the optimal value function \eqref{marg} with $\ph(x,y):=g(y)$. Thus Corollary~\ref{remark:FundThm} yields
\begin{align*}
\partial(g\circ F)(\ox)=A^*\partial g(A\ox)
\end{align*}
under the qualification conditions in (i) or (ii) of Theorem~\ref{sub chainrule}.}
\end{Remark}\vspace*{-0.12in}

The results for subdifferentiation of optimal value functions obtained above allow us to derive yet another chain rule involving generalized convex functions with values in {\em ordered} LCTV spaces; cf.\ \cite{lemarie}. Given an ordering convex cone $Y_+\subset Y$, define the {\em ordering relation} $\prec$ on $Y$ as follows: $y_1\prec y_2$  for $y_1,y_2\in Y$ if and only if $y_2-y_1\in Y_+$.\vspace*{-0.03in}

Recall some definitions.  Given a function $\ph\colon Y\to\oR$, we say that it is {\em $Y_+$-nondecreasing} if $\ph(y_1)\le\ph(y_2)$ for $y_1\prec y_2$. Given a mapping $f\colon X\to Y$, we say that $f$ is {\em $Y_+$-convex} if
\begin{equation*}
f\big(\lambda x+(1-\lambda)u\big)\prec\lambda f(x)+(1-\lambda)f(u)\;\;\mbox{\rm for all }\;x,u\in X\;\;\mbox{\rm and }\;\lambda\in(0,1).
\end{equation*}\vspace*{-0.25in}

It is easy to observe the following property.\vspace*{-0.2in}

\begin{Proposition}{\bf(subgradients of nondecreasing convex functions on ordered spaces).}\label{lem:Yplus} Let $\ph\colon Y\to\oR$ be $Y_+$-nondecreasing. Then for any $\oy\in\dom\ph$ and $y^*\in\partial\ph(\oy)$ we have $\la y^*,z\ra\ge 0$ whenever $z\in Y_+$.
\end{Proposition}\vspace*{-0.12in}
{\bf Proof.} It follows for any $z\in Y_+$ that $\oy-z\prec\oy$. Thus
\begin{equation*}
\la y^*,-z\ra=\la y^*,\oy-z-\oy\ra\le\ph(\oy-z)-\ph(\oy)\leq 0,
\end{equation*}
which clearly implies that $\la y^*,z\ra\ge 0$.$\h$

Now we are ready to derive the aforementioned chain rule in ordered spaces from Theorem~\ref{AnYenThm42}. For simplicity we present only the result corresponding assumption (i) therein.\vspace*{-0.12in}

\begin{Theorem}{\bf(chain rule for nondecreasing compositions in ordered spaces).}\label{cor:Lemarie} Let $f\colon X\to Y$ be $Y_+$-convex, and let  $\ph\colon Y\to\oR$ be convex and $Y_+$-nondecreasing. If there exists $x\in X$ such that $\ph$ is finite and continuous at some point $y\in Y$ with $f(x)\prec y$, then
\begin{equation*}
\partial(\ph\circ f)(\ox)=\bigcup_{y^*\in\partial\ph(f(\ox))}\partial(y^*\circ f)(\ox)\;\;\mbox{\rm whenever }\;\ox\in{\rm dom}(\ph\circ f).
\end{equation*}
\end{Theorem}\vspace*{-0.1in}
{\bf Proof.} Define $F\colon X\tto Y$ by $F(x):=\{y\in Y|\;f(x)\prec y\}$. It is easy to check that the $Y_+$-convexity $f$ yields that the graph of $F$ is convex in $X\times Y$. Since $\ph$ is $Y_+$-nondecreasing, we have the marginal function representation
\begin{equation*}
\mu(x):=\inf\big\{\ph(y)\big|\;y\in F(x)\big\}=(\ph\circ f)(x),\quad x\in X.
\end{equation*}
Then Theorem~\ref{AnYenThm42}(i) tells us that
\begin{align*}
\partial(\ph\circ f)(\ox)=\bigcup_{y^*\in\partial\ph(f(\ox))}D^*F\big(\ox,f(\ox)\big)(y^*).
\end{align*}
Then the claimed chain rule follows from the formula
\begin{equation}\label{ycoreq}
D^*F\big(\ox,f(\ox)\big)(y^*)=\partial(y^*\circ f)(\ox)\;\;\mbox{\rm whenever }\;y^*\in\partial\ph\big(f(\ox)\big).
\end{equation}
To verify first the inclusion ``$\subset$" in \eqref{ycoreq}, pick any $(y^*,x^*)\in\gph D^*F(\ox,f(\ox))$ and get that
\begin{equation*}
\la y^*,y\ra\ge\la y^*,f(\ox)\ra+\la x^*,x-\ox\ra\;\;\mbox{\rm for all }\;x\in X,\;f(x)\prec y.
\end{equation*}
Further, fix $h\in X$ and select $x:=\ox+h$, $y:=f(x)$. Since $f(x)\prec y$, it follows that
\begin{equation*}
\la y^*,f(\ox+h)\ra\ge\la y^*,f(\ox)\ra+\la x^*,h\ra,
\end{equation*}
which shows that $x^*\in\partial(y^*\circ f)(\ox)$, and hence the inclusion ``$\subset$" in \eqref{ycoreq} holds.\vspace*{-0.03in}

To verify the opposite inclusion therein, pick any $x^*\in\partial(y^*\circ f)(\ox)$ with $y^*\in\partial\ph(f(\ox))$ for which we have by the subgradient definition that
\begin{equation*}
\la y^*,f(\ox+h)-f(\ox)\ra\ge\la x^*,h\ra\;\;\mbox{\rm whenever }\;h\in X.
\end{equation*}
Taking any $x\in X$ and $y\in Y$ with $f(x)\prec y$, denote $h:=x-\ox$ and hence get $f(\ox+h)=f(x)$. Then Proposition~\ref{lem:Yplus} tells us that $\la y^*,y\ra\ge\la y^*,f(x)\ra$. This shows therefore that
\begin{equation*}
\la y^*,y-f(\ox)\ra\ge\la y^*,f(x)-f(\ox)\ra\ge\la x^*,x-\ox\ra
\end{equation*}
and so means by definition that $(x^*,-y^*)\in N((\ox,f(\ox));\gph F)$. Thus $x^*\in D^*F(\ox,f(\ox))(y^*)$, which justifies the inclusion ``$\supset$" in \eqref{ycoreq} and thus completes the proof of the theorem. $\h$\vspace*{-0.2in}

\section{Coderivative Calculus for Set-Valued Mappings}
\setcounter{equation}{0}\vspace*{-0.1in}

The final section of the paper is devoted to deriving new results of the coderivative calculus by using the geometric approach mainly based on the normal cone intersection rule  from Theorem~\ref{nir}. Note that the results obtained below are {\em essentially different} in various aspects from the coderivative calculus rules derived in \cite{m-book1,RockWets-VA} in the general nonconvex case as well as from their convex-graph specifications.\vspace*{-0.03in}

Let us start with the {\em sum rule} for coderivatives under several qualification conditions induced by those in Theorem~\ref{nir}(i--iii). For simplicity of formulations we confine ourselves in what follows to the coderivative implementation of the corresponding conditions formulated in (ii) and (iii) of Theorem~\ref{nir}.\vspace*{-0.03in}

Given two set-valued mappings $F_1,F_2\colon X\tto Y$ between LCTV spaces, consider their {\em sum}
\begin{equation*}
(F_1+F_2)(x)=F_1(x)+F_2(x):=\big\{y\in Y\big|\;\exists\,y_i\in F_i(x),\;i=1,2,\;\mbox{\rm such that }y=y_1+y_2\big\},\;x\in X.
\end{equation*}
It is easy to check that if $F_1$ and $F_2$ have convex graphs, then $F_1+F_2$ enjoys the same property, and also that ${\rm dom}(F_1+F_2)=\dom F_1\cap\dom F_2$. For $(\ox,\oy)\in{\rm gph}(F_1+F_2)$, define
\begin{equation*}
S(\ox,\oy):=\big\{(\oy_1,\oy_2)\in Y\times Y\big|\;\oy=\oy_1+\oy_2\;\;\mbox{\rm with }\;\oy_i\in F_i(\ox),\;i=1,2\big\}.
\end{equation*}\vspace*{-0.4in}

\begin{Theorem}{\bf(coderivative sum rule).}\label{thm:sum} Let $F_1,F_2\colon X\tto Y$ be set-valued mappings with convex graphs. Assume that one of the following conditions {\rm(i)}, {\rm(ii)} is satisfied:\\[1ex]
{\bf(i)} either ${\rm int}(\gph F_1)\ne\emp$ and ${\rm int}(\dom F_1)\cap\dom F_2\ne\emp$,\\
{\bf(ii)} or $X$ is Banach, the mappings $F_1$ and $F_2$ have closed graphs, and the set $\R^+(\dom F_1-\dom F_2)$ is a closed subspace of $X$.\\[1ex]
Then for any $(\ox,\oy)\in{\rm gph}(F_1+F_2)$ and $y^*\in Y^*$ we have
\begin{equation}\label{cod-sum}
D^*(F_1+F_2)(\ox,\oy)(y^*)=D^*F_1(\ox,\oy_1)(y^*)+D^*F_2(\ox,\oy_2)(y^*)\;\mbox{\rm if }(\oy_1,\oy_2)\in S(\ox,\oy).
\end{equation}
\end{Theorem}\vspace*{-0.12in}
{\bf Proof.} Let $y^*\in Y^*$, $(\ox,\oy)\in{\rm gph}(F_1+F_2)$, and $(\oy_1,\oy_2)\in S(\ox,\oy)$ be fixed for the entire proof. To verify first the inclusion ``$\subset$" in \eqref{cod-sum}, pick any $x^*\in D^*(F_1+F_2)(\ox,\oy)(y^*)$ and get by  the coderivative definition \eqref{cod} that
$(x^*,-y^*)\in N\big((\ox,\oy);{\rm gph}(F_1+F_2)\big)$. Define the convex subsets $\Omega_1$ and $\Omega_2$ of $X\times Y\times Y$ as follows:
\begin{equation}\label{cod-sum1}
\Omega_1:=\big\{(x,y_1,y_2)|\;y_1\in F_1(x)\big\},\;\Omega_2:=\big\{(x,y_1,y_2)\big|\;y_2\in F_2(x)\big\}.
\end{equation}
It is easy to check that $(x^*,-y^*,-y^*)\in N((\ox,\oy_1,\oy_2);\Omega_1\cap\Omega_2)$. By construction we get
\begin{equation}\label{graphrelation}
\Omega_1-\Omega_2=\big(\dom F_1-\dom F_2\big)\times Y\times Y.
\end{equation}
Imposing now the assumptions in (i) of the theorem yields $(\i\Omega_1)\cap\Omega_2\ne\emp$ for the sets $\O_i$ in \eqref{cod-sum}. Indeed, we have ${\rm int}(\Omega_1-\Omega_2)={\rm int}\big(\dom F_1-\dom F_2)\times Y\times Y.$ Since ${\rm int}(\dom F_1)\cap\dom F_2\ne\emp$, it follows that $0\in{\rm int}(\dom F_1-\dom F_2)$, and thus $0\in{\rm int}(\Omega_1-\Omega_2)$. Observe that ${\rm int}(\Omega_1)\ne\emp$ under (i). This yields ${\rm int}(\Omega_1)\cap \Omega_2\ne\emp$, since otherwise leads us to a contradiction by the separation theorem. Furthermore, the equality in \ref{graphrelation} implies that $\mathbb{R}^{+}(\Omega_1-\Omega_2)=\mathbb{R}^{+}(\dom F_1-\dom F_2)\times Y\times Y$, and so $\R^+(\Omega_1-\Omega_2)$ is a closed subspace under the validity of (ii). Applying now Theorem~\ref{nir} to the set $\O_i$ in \eqref{cod-sum1} gives us
\begin{equation*}
(x^*,-y^*,-y^*)\in N\big((\ox,\oy_1,\oy_2);\Omega_1\big)+N\big((\ox,\oy_1,\oy_2);\Omega_2\big),
\end{equation*}
which means that $(x^*,-y^*,-y^*)=(x^*_1,-y^*,0)+(x^*_2,0,-y^*)$ for some $x^*_1,x^*_2$ satisfying
\begin{equation*}
(x^*_1,-y^*)\in N\big((\ox,\oy_1);\gph F_1\big)\quad\text{and}\quad(x^*_2,-y^*)\in N\big((\ox,\oy_2);\gph F_2\big).
\end{equation*}
Thus $x^*=x^*_1+x^*_2\in D^*F_1(\ox,\oy_1)(y^*)+D^*F_2(\ox,\oy_2)(y^*)$, and ``$\subset$" in \eqref{cod-sum} is justified.\vspace*{-0.03in}

To verify the opposite inclusion in \eqref{cod-sum}, pick $x^*\in D^*F_1(\ox,\oy_1)(y^*)+D^*F_2(\ox,\oy_2)(y^*)$ and find $x^*_1\in D^*F_1(\ox,\oy_1)(y^*)$ and $x^*_2\in D^*F_2(\ox,\oy_2)(y^*)$ with $x^*=x^*_1+x^*_2$. This tells us by the coderivative and normal cone definitions that
\begin{equation*}
\la x^*_1,x_1-\ox\ra+\la x^*_2,x_2-\ox\ra-\la y^*,y_1-\oy_1\ra-\la y^*,y_2-\oy_2\ra\le 0\;\;\mbox{\rm for all }\;y_i\in F_i(x_i),\;i=1,2.
\end{equation*}
Setting $x_1=x_2=:x$ gives us the relationships
\begin{equation*}
\la w,x-\ox\ra-\la y^*,y_1-\oy_1\ra-\la y^*,y_2-\oy_2\ra\le 0\;\;\mbox{\rm whenever }\;y_1+y_2\in(F_1+F_2)(x).
\end{equation*}
Denoting finally $y:=y_1+y_2$, we arrive at the inequality
\begin{equation*}
\la x^*,x-\ox\ra-\la y^*,y-\oy\ra\le 0\;\;\mbox{\rm for all }\;(x,y)\in{\rm gph}(F_1+F_2),
\end{equation*}
which tells us that $x^*\in D^*(F_1+F_2)(\ox,\oy)(y^*)$ and thus completes the proof. $\h$

Our next goal is to derive geometrically a precise {\em coderivative chain rule} for compositions of general set-valued mappings between LCTV spaces. Given $F\colon X\tto Y$ and $G\colon Y\tto Z$, recall that their {\em composition} $(G\circ F)\colon X\tto Z$ is defined by
\begin{equation*}
(G\circ F)(x)=\bigcup_{y\in F(x)}G(y):=\big\{z\in Z\big|\;\exists y\in F(x)\;\;\mbox{\rm with }\;z\in G(y)\big\},\quad x\in X.
\end{equation*}
It is easy to see that if $F$ and $G$ have convex graphs, then $G\circ F$ also has this property. In the following theorem we use the notation
\begin{equation*}
T(x,z):=F(x)\cap G^{-1}(z)\;\;\mbox{\rm and }\;\rge F:=\cup_{x\in X}F(x).
\end{equation*}\vspace*{-0.4in}

\begin{Theorem}{\bf(chain rule for coderivatives).}\label{thm:comp} Let $F\colon X\tto Y$ and $G\colon Y\tto Z$ be convex-graph mappings. Suppose that one of the following groups of assumptions is satisfied:\\[1ex]
{\bf(i)} Either ${\rm int}(\gph F)\ne\emp$ and ${\rm int}(\rge F)\cap\dom G\ne\emp$, or ${\rm int}(\gph F)\ne\emp$ and the condition $\mbox{\rm rge}\,F\cap({\rm int}(\dom G))\ne\emp$ holds.\\
{\bf(ii)} The spaces $X,Y,Z$ are Banach, the set-valued mappings $F$ and $G$ are closed-graph, and the set $\R^+(\rge F-\dom G)$ is a closed subspace of $Y$.\\[1ex]
Then for all $(\ox,\oz)\in\gph(G\circ F)$ and for all $z^*\in Z^*$ we have the coderivative chain rule
\begin{equation}\label{cod-chain}
D^*(G\circ F)(\ox,\oz)(z^*)=\big(D^*F(\ox,\oy)\circ D^*G(\oy,\oz)\big)(z^*)\;\;\mbox{\rm whenever }\;\oy\in T(\ox,\oz).
\end{equation}
\end{Theorem}\vspace*{-0.15in}
{\bf Proof.} Let $(\ox,\oz)\in{\rm gph}(G\circ F)$, $z^*\in Z^*$, and $\oy\in T(\ox,\oz)$ be fixed for the entire proof. We first verify the inclusion ``$\subset$" in \eqref{cod-chain}. Pick $x^*\in D^*(G\circ F)(\ox,\oz)(z^*)$ and get
\begin{align*}
\la x^*,x-\ox\ra-\la z^*,z-\oz\ra\le 0\;\;\mbox{\rm for all }\;(x,z)\in{\rm gph}(G\circ F).
\end{align*}
Define the following subsets of $X\times Y\times Z$:
\begin{equation}\label{cod-chain1}
\Omega_1:=\big\{(x,y,z)\big|\;(x,y)\in\gph F\big\},\quad\Omega_2:=\big\{(x,y,z)\big|\;(y,z)\in\gph G\big\}.
\end{equation}
It is easy to see that $(x^*,0,-z^*)\in N((\ox,\oy,\oz);\Omega_1\cap\Omega_2)$ and that
\begin{equation*}\label{graphrelation1}
\Omega_1-\Omega_2=X\times(\rge F-\dom G)\times Z.
\end{equation*}
Similar to the proof of Theorem \ref{thm:sum} we conclude that the assumptions in (i) imply that $(\i\Omega_1)\cap\Omega_2\ne\emp$. It also follows that $\mathbb{R}^{+}(\Omega_1-\Omega_2)=X\times\mathbb{R}^{+}(\rge F-\dom G)\times Z$, and so $\Omega_1-\Omega_2$ is a closed subspace of $X\times Y\times Z$ under (ii).
Furthermore, we easily deduce from the structures of $\O_1,\O_2$ in \eqref{cod-chain1} the normal cone representations:
\begin{align*}
N\big((\ox,\oy,\oz);\Omega_1\big)=N\big((\ox,\oy);\gph F\big)\times\{0\},\quad N\big((\ox,\oy,\oz);\Omega_2\big)=\{0\}\times N\big((\oy,\oz);\gph G\big).
\end{align*}
It gives us $(x^*,0,-z^*)=(x^*,-y^*,0)+(0,y^*,-z^*)$ with $(x^*,-y^*)\in N((\ox,\oy);\gph F)$ and $(y^*,-z^*)\in N((\oy,\oz);\gph G)$.  This means that $y^*\in D^*G(\oy,\oz)z^*)$ and $x^*\in D^*F(\ox,\oy)(y^*)$ and thus yields $x^*\in(D^*F(\ox,\oy)\circ D^*G(\oy,\oz))(z^*)$, i.e., justifies the inclusion ``$\subset$" in \eqref{cod-chain}.\vspace*{-0.03in}

To verify the opposite inclusion therein, pick $x^*\in(D^*F(\ox,\oy)\circ D^*G(\oy,\oz))(z^*)$ and find $y^*\in D^*G(\oy,\oz)(z^*)$ such that $x^*\in D^*F(\ox,\oy)(y^*)$. It tells us that
\begin{align*}
\la x^*,x-\ox\ra-\la y^*,y_1-\oy\ra\le 0\;\text{ whenever }\;y_1\in F(x),\\
\la y^*,y_2-\oy\ra-\la z^*,z-\oz\ra\le 0\;\text{ whenever }\;z\in G(y_2).
\end{align*}
Summing up these two inequalities and putting $y_1=y_2=:y$, we get
\begin{align*}
\la x^*,x-\ox\ra-\la z^*,z-\oz\ra\le 0 \;\text{ for all }\;z\in(G\circ F)(x),
\end{align*}
which justifies $x^*\in D^*(G\circ F)(\ox,\oz)(z^*)$ and completes the proof of the theorem. $\h$

It follows from the proof of Theorem~\ref{thm:comp} that in assumption (i) therein the condition ${\rm int}(\gph F)\ne\emp$ can be replaced by ${\rm int}(\gph G)\ne\emp$.\vspace*{-0.03in}

The final result of this section provides a useful rule for representing the coderivative of intersections of set-valued mappings. Again we derived it geometrically from the normal cone intersection rule while concentrating for simplicity on applications of the easily verifiable conditions in (ii) and (iii) of Theorem~\ref{nir}.  Given two set-valued mappings
$F_1,F_2\colon X\tto Y$ between LCTV spaces, recall that their {\em intersection} $(F_1\cap F_2)\colon X\tto Y$ is defined by
\begin{equation*}
(F_1\cap F_2)(x):=F_1(x)\cap F_2(x),\quad x\in X.
\end{equation*}
It is easy to see that ${\rm gph}(F_1\cap F_2)=({\rm gph}\,F_1)\cap({\rm gph}F_2)$, and so the convexity of both sets $F_1$ and $F_2$ yields the convexity of ${\rm gph}(F_1\cap F_2)$.  This allows us to derive the following intersection rule for coderivatives.\vspace*{-0.12in}

\begin{Theorem}{\bf(coderivative intersection rule).}\label{thm:int} Let $F_1,F_2\colon X\tto Y$ be set-valued mappings with convex graphs. Suppose that one of the conditions {\rm(i)}, {\rm(ii)} is satisfied:

{\bf(i)} ${\rm int}(\gph F_1)\cap(\gph F_2)\ne\emp$.\\
{\bf(ii)} The spaces $X$ and $Y$ are Banach, the mappings $F_1$ and $F_2$ are closed-graph, and the set $\R^+(\gph F_1-\gph F_2)$ is a closed subspace of $X\times Y$.\\[1ex]
Then for any $\oy\in(F_1\cap F_2)(\ox)$ and any $y^*\in Y^*$ we have
\begin{equation}\label{cod-inter}
D^*(F_1\cap F_2)(\ox,\oy)(y^*)=\bigcup_{y^*_1+y^*_2=y^*}\big(D^*F_1(\ox,\oy)(y^*_1)+D^*F_2(\ox,\oy)(y^*_2)\big).
\end{equation}
\end{Theorem}\vspace*{-0.1in}
{\bf Proof.} First we verify the inclusion ``$\subset$" in \eqref{cod-inter}. For every $\oy\in(F_1\cap F_2)(\ox)$, $y^*\in Y^*$, and $x^*\in D^*(F_1\cap F_2)(\ox,\oy)(y^*)$ it follows that
\begin{equation*}
(x^*,-y^*)\in N\big((\ox,\oy);{\rm gph}(F_1\cap F_2)\big)=N\big((\ox,\oy);(\gph F_1)\cap(\gph F_2)\big).
\end{equation*}
Then applying Theorem~\ref{nir} under the corresponding assumptions tells us that
\begin{align*}
(x^*,-y^*)\in N\big((\ox,\oy);{\rm gph}(F_1\cap F_2)\big)=N\big((\ox,\oy);\gph F_1\big)+N\big((\ox,\oy);\gph F_2\big).
\end{align*}
Thus $(x^*,-y^*)=(x^*_1,-y^*_1)+(x^*_2,-y^*_2)$ with some $(x^*_i,-y^*_i)\in N((\ox,\oy);\gph F_i))$ for $i=1,2$. Therefore $x^*\in D^*F_1(\ox,\oy)(y^*_1)+D^*F_2(\ox,\oy)(y^*_2)$ and $y^*=y^*_1+y^*_2$, which justify the claimed inclusion ``$\subset$" in the coderivative representation \eqref{cod-inter}.\vspace*{-0.03in}

To verify the opposite inclusion in \eqref{cod-inter}, take $y^*_1,y^*_2\in Y^*$ with $y^*_1+y^*_2=y^*$. Picking now $x^*\in D^*F_1(\ox,\oy)(y^*_1)+D^*F_2(\ox,\oy)(y^*_2)$, we get $x^*=x^*_1+x^*_2$ with some $x^*_1\in D^*F_1(\ox,\oy)(y^*_1)$ and $x^*_2\in D^*F_2(\ox,\oy)(y^*_2)$. This shows that
\begin{align*}
(x^*,-y^*)=(x^*_1,-y^*_1)+(x^*_2,-y^*_2)&\in N\big((\ox,\oy);\gph F_1\big)+N\big((\ox,\oy);\gph F_2\big)\\
&=N\big((\ox,\oy);{\rm gph}(F_1\cap F_2)\big),
\end{align*}
and thus $x^*\in D^*(F_1\cap F_2)(\ox,\oy)(y^*)$, which completes the proof of the theorem. $\h$\vspace*{0.05in}

{\bf Acknowledgements.} The authors are grateful to both anonymous referees for their helpful comments that allowed us to improve the original presentation. The authors are also grateful to Bingwu Wang for helpful discussions on the material presented in this paper. 
\end{document}